\documentclass{article}

\usepackage{arxiv}

\usepackage[utf8]{inputenc} 
\usepackage[T1]{fontenc}    
\usepackage{hyperref}       
\usepackage{url}            
\usepackage{booktabs}       
\usepackage{amsfonts}       
\usepackage{nicefrac}       
\usepackage{microtype}      
\usepackage{lipsum}		
\usepackage{graphicx}
\usepackage{natbib}
\usepackage{doi}
\usepackage{amsmath}
\usepackage{algorithm}
\usepackage{algpseudocode}
\usepackage{amssymb}

\title{Contribution of expert aggregation to temperature prediction part I.}

\author{ {Leo ~Pfitzner} \\
	Meteo France\\
	\texttt{leo.pfitzner@meteo.fr} \\
	\And
	{Olivier~Wintenberger} \\
	Sorbonne University\\
	\texttt{olivier.wintenberger@sorbonne-universite.fr} \\
    \And
	{Olivier~Mestre} \\
	Meteo France\\
	\texttt{olivier.mestre@meteo.fr} \\
    \And
	{Marion~Riverain} \\
	Meteo France\\
}

\renewcommand{\shorttitle}{\textit{arXiv} Template}

\hypersetup{
pdftitle={A template for the arxiv style},
pdfsubject={q-bio.NC, q-bio.QM},
pdfauthor={David S.~Hippocampus, Elias D.~Striatum},
pdfkeywords={First keyword, Second keyword, More},
}

\begin{document}
\maketitle

\begin{abstract}
	Many Numerical Weather Prediction models and their associated Post-Processed Models are available. Combining all of these predictions in an optimal way is however not straightforward. This can be achieved thanks to Expert Aggregation (EA) which has many advantages, such as being online, being adaptive to model changes and having theoretical guarantees. In this paper, we propose a method for making deterministic temperature predictions with EA. We used Exponentially Weighted Average, MLprod and MLpol and Bernstein Online Aggregation. Hence, we combine and outperform the forecasts of the raw and post-processed Integrated Forecasting System (IFS), forecasts of Application of Research to Operations at Mesoscale (AROME), Action de Recherche Petite Echelle Grande Echelle (ARPEGE) and quantiles of the post processed Prévision d’Ensemble ARPEGE (PEARP). We also compare the different EA strategies in various settings and show that they outperform the National Blend of Models. Finally, we discuss certain limitations.
\end{abstract}


\section{Introduction}
\hspace{0.5cm}In order to predict weather and temperature, many Numerical Weather Prediction (NWP) models are available \citep{seity_arome-france_2011,barros_ifs_1995,deque_arpegeifs_1994}. Additionally, so called ensemble NWP models have been introduced in order to deal with the uncertainty of the forecast \citep{leutbecher_ensemble_2008}.

All these NWP models suffer however from systematic errors. The ensemble models also typically exhibit under-dispersion \citep{hamill_verification_1997}. This is why Post-Processed Models (PPM) are more and more used in order to reduce systematic errors thanks to machine learning techniques \citep{vannitsem_statistical_2021}.

But these machine learning techniques struggle to deal with changing environments, especially with changes in the NWP models themselves. For this reason they often need to be retrained when the raw models are modified because they are "static" and not "online". By "online" model, we mean a model that is updated at each prediction, contrary to usual static machine learning techniques.

Confronted with all of these - raw NWP models and post-processed NWP models - the natural question arises whether it is possible \i) to combine/aggregate them in an optimal way so as to improve the predictions compared to all the models that we use, \i\i) to do this using an online method in order to avoid prior training and retraining when the environment changes.

In this paper, we do not consider aggregation methods that are not online, such as multiple regression \citep{krishnamurti_improved_1999,krishnamurti_multimodel_2000,kharin_climate_2002}. Additionally, we do not address methods designed for probabilistic forecasting, such as Bayesian model averaging \citep{raftery_using_2005,wilson_calibrated_2007}, since we consider deterministic forecasts only. Finally, seamless forecasts such as the operational IMPROVER framework \citep{roberts_improver_2023} are beyond the scope of this work. 

As to online aggregations built to improve deterministic forecasts, first investigations simply took the mean - the uniform aggregation - of several NWP models as in \citet{ebert_ability_2001} for 24h and 48h precipitations. In \citet{mylne_multi-model_2002}, they computed the mean of all the members of two ensemble NWP models - for mean sea level pressure, 500 hPa geopotential height, 850 hPa temperature and daily precipitation accumulations.

A more sophisticated online deterministic aggregation method consists in reducing the bias of the models before averaging them with online methods. The bias correction was performed by \citet{delle_monache_ozone_2006} for ozone predictions using Kalman filters, which were subsequently averaged. In the Operational Consensus Forecasting (OCF) \citep{woodcock_operational_2005,engel_gridded_2012} they used the best easy systematic mean statistic to reduce the bias followed by a mean absolute error-weighted average. For the National Blend of Models (NBM) of \citet{craven_national_2020}, the authors used the "decaying averaging method" \citep{cui_bias_2012} to reduce the bias of temperature predictions among other parameters like wind speed. Then, they applied a decaying averaging method to the mean absolute error to derive a weighted mean.


NBM \citep{craven_national_2020} is an operational online aggregation method of the National Weather Service in the USA. It will be our main benchmark.

In this study, we investigate how to aggregate the different models using the Expert Aggregation (EA) approach described in \citet{cesa-bianchi_prediction_2006}. This framework provides simple online algorithms with theoretical guarantees to combine predictions from arbitrary forecasting sources, referred to as "experts".

EA has been introduced by \citet{vovk_aggregating_1990} and \citet{littlestone_weighted_1994} with the Exponentially Weighted Average (EWA) in a very general framework, where no assumptions are made on the observations nor on the experts which are seen as black boxes. Recently, more subtle EAs have been developed by \citet{gaillard_second-order_2014} (MLprod, and MLpol), \citet{van_erven_metagrad_2021} (MetaGrad), \citet{koolen2015second} (Squint), \citet{wintenberger_optimal_2017} and \citet{wintenberger_stochastic_2024} (Bernstein Online Aggregation (BOA)) among others. 

EA has already proved to be useful in several different fields in practice. For example, for ozone predictions \citep{mallet_ozone_2009}, electricity consumption \citep{devaine_forecasting_2013}, finance \citep{amat_fundamentals_2018,remlinger_expert_2023}, oil and gas \citep{deswarte_sequential_2019}, probabilistic photovoltaic forecasts \citep{thorey_ensemble_2018}, wind farm forecast \citep{van_der_meer_crps-based_2024} probabilistic wind forecasts \citep{zamo_sequential_2021} and probabilistic synthetic temperature prediction \citep{thorey_online_2017}.

 At Meteo France, BOA has been used operationally for several years with a 60 days sliding window to predict 2 meters above the ground temperature. The workflow is as follows: first, the NWP models are run; then they are post-processed; finally, the expert aggregation method combines both the raw models and the PPMs. 

In this study, we combine the raw and post-processed Integrated Forecasting System (IFS), forecasts of Application of Research to Operations at Mesoscale (AROME), Action de Recherche Petite Echelle Grande Echelle (ARPEGE) and quantiles of the post processed Prévision d’Ensemble ARPEGE (PEARP). We use raw models but also PPMs as they perform better on average, and we also want to improve them further. Hence, we post-process models that have already been post-processed, which represents a significant challenge in terms of making further improvements. To do so, we will use EWA, BOA, MLprod and MLpol and compare them to NBM and to the experts and their convex combinations in the spirit of EA.

Our main result is that one can outperform NBM and all the individual experts including the PPMs by using EA (see Table \ref{rmse_ewa_boa_MLpol_with_grad_with_pearp}). Furthermore, to the best of our knowledge, it is the first time that EA for deterministic temperature forecasting is done with raw NWP models, PPMs and quantiles of a post processed ensemble NWP model. We also think that it is the first time that aggregation strategies which are adapted for worst cases (EWA) are compared with EAs which are more adaptive (BOA, MLpol and MLprod), and should be more suited to temperature forecasting.

The paper is organized as follows. In Section \ref{section_theoretical_framework}, we present the theoretical framework, in particular we present the algorithms of EWA, BOA, MLprod and MLpol along with their theoretical guarantees. In the third Section, we present the data and experts that we used, how we evaluated the forecasts and we present our benchmarks. In the fourth Section, we first examine the influence of the set of experts on EA. Then we show how EA improves temperature predictions, especially compared to the best expert, the uniform aggregation and the best convex combination of experts and NBM. Finally, we will point out some limitations and see if it is possible to handle these limitations with sliding windows. Section \ref{section_conclusion} contains the conclusion which is followed by the Appendix.

\section{Theoretical framework}
\label{section_theoretical_framework}
\subsection{Expert Aggregation}
\label{subsection:expert_agregation}

For a summary of the literature on sequential prediction we refer to the seminal monograph \citet{cesa-bianchi_prediction_2006} (in particular to Chapter 2).

The goal of EA is to predict a sequence of observations $y_1,\ldots,y_T$ within an outcome set $\mathcal{Y}\subset\mathbb{R}$ for a number $T$ of iterations. In our case the sequence to be predicted is the temperature at 2 meters above ground for a given station and lead time.

In order to perform this prediction, at each iteration $t=1,\ldots,T$, a number $N$ of what we call experts which can be seen as black boxes, make a prediction $x_{1,t},\ldots,x_{N,t}$ within a decision set $\widehat{\mathcal{Y}}\subset{\mathbb{R}}$. Concerning our application, these experts are the raw and post-processed NWP models (outlined in subsection \ref{section_data_evaluation}.\ref{subsection_data_experts}), for the specific lead time and station considered.

This enables the aggregation at $t\geqslant1$ to make its own prediction $\widehat{y}_t \in \widehat{\mathcal{Y}}$, by providing a weight vector $\mathbf{w}_t=(w_{1,t},\ldots,w_{N,t})\in \mathbb{R}^N$ and putting $\widehat{y}_t=\sum_{i=1}^N w_{i,t} x_{i,t}$. The weight vector $\mathbf{w}_t=(w_{1,t},\ldots,w_{N,t})\in \mathbb{R}^N$ only depends on present and past information. $\mathbf{w}_t$ will especially depend on $y_s$ with $s<t$ and $x_{i,s}$with $s\leqslant t$, $i=1,\ldots, N$. In our study we will focus on convex aggregation strategies, such that $\mathbf{w}_t\in \mathcal{X} \subset \mathbb{R}^N$ where $\mathcal{X}$ is the set of convex weight vectors $\mathbf{w}_t$ such that $\sum_{i=1}^N w_{i,t} = 1$ and $w_{i,t} \geqslant 0$ for all $i=1,\ldots,N$. Using convex aggregation strategies implies two important aspects. The predictions are easy to interpret for human forecasters and the predictions are between the maximum and the minimum of the experts’ predictions.

After the aggregation's prediction, the true observation $y_t$ is revealed. The aggregation incurs the loss $\ell(\widehat{y}_t,y_t)=\ell(\sum_{i=1}^N w_{i,t}x_{i,t},y_t)$ and each expert $i$ incurs the loss $\ell(x_{i,t},y_t)$ with $\ell:\widehat{\mathcal{Y}} \times\mathcal{Y} \rightarrow \mathbb{R}$ denoting the loss function. We set $\ell_t:\mathcal{X} \rightarrow \mathbb{R}$ by $\ell_t(\mathbf{q})=\ell(\sum_{i=1}^N q_ix_{i,t},y_t)$, where the subscript $t$ of $\ell_t$ refers to the expert predictions $x_{1,t},\ldots,x_{N,t}$ and the observation $y_t$. Hence, $\ell_t(\mathbf{w}_t)=\ell(\widehat{y}_t,y_t)$ and $\ell_t(\mathbf{\delta}_i)=\ell(x_{i,t},y_t)$, where $\mathbf{\delta}_i\in \mathbb{R}^N$ is the Dirac mass i.e. the weight vector with 1 on the $i^{\text{th}}$ coordinate and zero elsewhere.

Once the losses are known, the weight vector is updated for the next iteration, increasing (or decreasing) the weights of the experts which performed well (or badly) at $t$. Algorithm \ref{expert_agregation_algorithm} summarizes EA.

\begin{algorithm}
\caption{Expert aggregation (chapter 2 of \citet{cesa-bianchi_prediction_2006})}
\label{expert_agregation_algorithm}
\begin{algorithmic}
    \State N experts.\\
    \For{$t \geqslant 1$}

    \begin{itemize}
        \item Expert $i$ makes the prediction $x_{i,t}$, $i=1,\ldots,N$.
        \item The aggregation makes the prediction $\widehat{y}_{t}=\sum_{i=1}^N w_{i,t} x_{i,t}$.
        \item The environment reveals the true outcome $y_{t}$.
        \item The aggregation incurs the loss $\ell_{t}(\mathbf{w}_t)$ and expert $i$ the loss $\ell_{t}(\mathbf{\delta}_i)$, $i=1,\ldots,N$.
        \item The aggregation updates $\mathbf{w}_{t}$ to $\mathbf{w}_{t+1}$.
    \end{itemize}
    
    \EndFor
\end{algorithmic}
\end{algorithm}

At each iteration $t$, the cumulative loss of the aggregation is defined to be $\mathcal{L}_{t}=\sum_{s=1}^t \ell_{s}(\mathbf{w}_s)$ and the one of a convex combination $\mathbf{q}\in\mathcal{X}$ of experts is $\mathcal{L}_{t}(\mathbf{q})=\sum_{s=1}^t \ell_{s}(\mathbf{q})$.
The cumulative loss of expert $i$ is $\mathcal{L}_{t}(\mathbf{\delta}_i)=\sum_{s=1}^t \ell_{s}(\mathbf{\delta}_i)$.

Throughout our applications, we will use only the square loss $\ell(x,y)=(x-y)^2$ for $x \in \widehat{\mathcal{Y}}$ and $y \in \mathcal{Y}$ where $\widehat{\mathcal{Y}}$ and $\mathcal{Y}$ are bounded intervals of $\mathbb{R}$. Note in passing that the square loss function is positive, convex and differentiable in the first variable.

The initial goal of the aggregation (as of any forecaster) is to have the cumulative loss as small as possible. In the worst case, the adversarial\footnote{The adversarial setting comes from game theory. It is a game where a forecaster (our aggregation) plays against someone or something (in our case, it would be the atmosphere) which wants the aggregation to be as bad as possible. To achieve this, the adversary can choose the observation/the losses of the experts, in order to make the forecasters’ loss as large as possible. This framework is too pessimistic for temperature prediction. Indeed, the player we are "playing against" - the atmosphere - does not care about our predictions.} case, when all the experts may have large losses, the loss of the aggregation can be arbitrarily high. Therefore, instead of minimizing its cumulative loss, the aggregation will try to minimize what is called its regret $\mathcal{R}_{T}$.

The regret $\mathcal{R}_{T}$ is the difference between the aggregation's cumulative loss and the cumulative loss of a reference forecast. In our case, depending on the aggregation strategy, this reference forecaster can be any element of the set of experts or, more generally, of the set of convex combinations of experts. The best strategy of a set of forecaster is called "oracle", since the latter is revealed at the end of the forecast sequence. We will use the notion of best expert (or best convex combination) on average, which will be made precise in Subsection \ref{section_data_evaluation}.\ref{benchmarks}.

This leads to the definition of the regret of the aggregation strategy against a convex combination $\mathbf{q} \in \mathcal{X}$ of experts,

\begin{equation}
 \mathcal{R}_{T}(\mathbf{q})= \mathcal{L}_T-\mathcal{L}_T(\mathbf{q})= \sum_{t=1}^T \ell_{t}(\textbf{w}_t) - \ell_{t}(\mathbf{q}).\label{eq def-regret}
\end{equation}

In particular, the regret against expert $i$ is $\mathcal{R}_{T}(\mathbf{\delta}_i)=\mathcal{L}_T-\mathcal{L}_T(\mathbf{\delta}_i)$, the regret against the best expert is $\max_{i=1,\ldots,N}\mathcal{R}_{T}(\mathbf{\delta}_i)$ and the regret against the best combination of experts is $\max_{\mathbf{q}\in\mathcal{X}}\mathcal{R}_{T}(\mathbf{q})$.

The goal of an aggregation strategy is to have $\mathcal{R}_t$ as small as possible. Hence, a suitable property is when the regret is equal or smaller than a "little o"\footnote{Let $f,g:\mathbb{R}\to\mathbb{R}$ be two functions with $g>0$ on an interval of the form $]a;\infty[$ for a certain real $a$. We say that $f$ is a little $o$ of $g$ and write $f(t)=o(g(t))$ if $\lim_{t\rightarrow+\infty}f(t)/g(t)=0$} of T, $\mathcal{R}_{T}\leqslant o(T) $. Equivalently and more simply, this is the case when the average performance given by the mean regret $\frac{1}{T} \mathcal{R}_{T}$
converges to zero or even less.


In this case, without prior knowledge on the oracle, the aggregation will at least perform as well on average as the oracle, when $T$ increases.  

The condition $\mathcal{R}_{T}\leqslant o(T)$ is what we meant by theoretical guarantees in the introduction, see inequalities \eqref{eq:EWA_optimal_regret}, \eqref{eq:EWA_grad_optimal_regret}, \eqref{eq:regret_BOA_bound}, \eqref{eq:regret_BOA_grad_bound}, \eqref{eq:regret_MLprod_bound}, \eqref{eq:regret_MLprod_grad_bound}, \eqref{eq:regret_MLpol_bound} and \eqref{eq:regret_MLpol_grad_bound} below. The fact that it is an asymptotic property will be implicit in the sequel.

\subsection{Exponentially Weighted Average}
\label{EWA_part}
One of the simplest aggregation strategies is the Exponentially Weighted Average (EWA) \citep{vovk_aggregating_1990,littlestone_weighted_1994}. EWA is described in Algorithm \ref{EWA_aggregation}. It is an aggregation strategy built to guarantee optimal minimax\footnote{Minimax is a term from game theory. Here, it means that it guarantees that it will minimize the maximum loss; in other words, it is certain that we will be as good as possible in worst/adversarial cases.} regret rates against the experts in an adversarial environment as explained in Section 2.10 of \citet{cesa-bianchi_prediction_2006}.

\begin{algorithm}
\caption{EWA aggregation \citep{vovk_aggregating_1990}}
\label{EWA_aggregation}
\begin{algorithmic}
    \State Parameter: $\eta>0$
    \State Initialization: $\mathbf{w_{1}} \in \mathbb{R}^N $ the uniform weight vector.
    \For{$t \geqslant 1$}

    \begin{equation}
        w_{i,t+1} = \frac{w_{i,t}e^{-\eta \ell_{t}(\mathbf{\delta}_i)}}{\sum_{j=1}^N w_{j,t}e^{-\eta \ell_{t}(\mathbf{\delta}_j)}},i=1,\ldots,N
    \end{equation}
    \EndFor
\end{algorithmic}
\end{algorithm}

If the loss function is convex in its first argument, positive and bounded by $M \in \mathbb{R}$, then after $T \geqslant 1$ iterations, the regret $\mathcal{R}_{T}(\mathbf{\delta}_i)$ of EWA against expert $i$ with the the learning rate $\eta>0$ is bounded as follows (see Theorem 2.2 of \citet{cesa-bianchi_prediction_2006} , \citet{stoltz_agregation_2010}):

\begin{equation}
    \mathcal{R}_T^\text{EWA}(\mathbf{\delta}_i) \leqslant \frac{\ln(N)}{\eta} + \frac{\eta}{8} M^2 T,
\end{equation}

When competing against the experts, the theoretical optimal choice for the learning rate  $\eta$ for $T$ iterations, $N$ experts and a loss bounded by $M\in \mathbb{R}$ is $\eta^*=M^{-1} \sqrt{(8\ln(N))/T}$. In this case, the regret of EWA against expert $i=1,\ldots,N$ becomes: 

\begin{equation}
\label{eq:EWA_optimal_regret}
    \mathcal{R}_T^\text{EWA}(\mathbf{\delta}_i) \leqslant M \sqrt{ \frac{T}{2} \ln(N)} = o(T)
\end{equation}

In this type of aggregation strategy with a learning rate $\eta$, the choice of $\eta$ is crucial and has a very important impact on the behavior of the aggregation. The theoretically optimal choice of $\eta^*$ is not always best in practice. For example, inequality \eqref{eq:EWA_optimal_regret} holds for worst-case settings, so $\eta^*$ may be too conservative in non-adversarial environments.  
The larger the learning rate $\eta$ is, the faster the aggregation learns and changes the weights of the experts.

When $\eta$ becomes very large, $\eta \longrightarrow \infty$, then EWA becomes the "Follow The Leader" aggregation and puts all the weight on the last best expert and predicts in the same way as this expert. And when $\eta \longrightarrow 0$ then EWA becomes the uniform aggregation. Therefore, in practice as mentioned in \citet{koolen_learning_2014} one has to choose an $\eta$ which is large enough to fit the data well, but small enough to prevent overfitting.

In EA, a well known trick to compete against convex combinations of experts is the gradient trick presented in \citet{cesa-bianchi_prediction_2006} Section 2.5. The idea is to give large weights to experts pointing in the direction where the loss function decreases, i.e. pointing in the opposite direction of the gradient.

When the loss function is convex and differentiable
in its first argument, then the gradient trick consists in replacing the loss $\ell_t(\mathbf{q})$ by the (linear) loss $\nabla \ell_t(\mathbf{w}_t) \cdot \mathbf{q}$ where:

\begin{equation}
    \nabla \ell_t(\mathbf{w}_t) \cdot \mathbf{q} = \frac{ \partial \ell(\widehat{y}_t,y_t)}{\partial \widehat{y}_t} \sum_{i=1}^N q_i x_{i,t},
\end{equation}

with $\mathbf{q}=(q_1,...,q_N) \in \mathcal{X}$ a convex weight vector, $\nabla \ell_t(\mathbf{w}_t)$ the gradient of  $\ell_t$ at $\mathbf{w}_t$, $y_t$ the observation, $\widehat{y}_t$ the prediction of the EA at $t$, and $\cdot$ the inner product in $\mathbb{R}^N$. For the square loss $\ell(\widehat{y}_t,y_t)=(\widehat{y}_t-y_t)^2$, this becomes:
\begin{equation}
    \nabla \ell_t(\mathbf{w}_t) \cdot \mathbf{q}=2\left(\sum_{i=1}^N w_{i,t}x_{i,t}-y_t\right) \sum_{i=1}^N q_i x_{i,t}=2(\widehat{y}_t-y_t) \sum_{i=1}^N q_i x_{i,t}.
\end{equation}

The convexity and differentiability of the loss $\ell_t$ imply the inequality:

\begin{equation}
\ell_{t}(\textbf{w}_t) - \ell_{t}(\mathbf{q})
\le\nabla \ell_t(\mathbf{w}_t) \cdot (\mathbf{w}_t-\mathbf{q})
=\nabla \ell_t(\mathbf{w}_t) \cdot (\mathbf{w}_t) -\nabla \ell_t(\mathbf{w}_t) \cdot (\mathbf{q}),\label{eq diff des pertes}
\end{equation}

 which enables to control the generic term in the definition of the regret \eqref{eq def-regret}.

Note in passing that the gradient trick still works if the loss function is convex and only subdifferentiable\footnote{The subgradient generalizes the gradient to functions that are convex but not necessarily differentiable.
If $f: U \longrightarrow \mathbb{R}$ is a real-valued convex function defined on a convex open set $U$ in the Euclidean space $\mathbb{R}^{n}$, then $v\in U$ is a subgradient of $f$ at $u\in U$ if for any $x\in U$ one has  $f(u)-f(x)\leqslant v\cdot (u-x)$. A fundamental result in convex optimization states that if $f$ is as just described then such a subdifferential exists at each $u\in U$. In general, it is not unique but it is when $f$ is differentiable, see e.g.\citet{borwein_convex_2006} Subsection 3.1.} in the first argument (e.g.\ the absolute loss function $\mathbb{R}^2\ni(x,y)\mapsto|x-y|$). For more details compare e.g.\ \citet{devaine_forecasting_2013} Subsection 2.2 and \citet{stoltz_agregation_2010} Subsection 2.3. 

In the case of a convex loss function, and when the linearized loss $\nabla \ell_t(\mathbf{w}_t) \cdot \mathbf{q} $ is bounded in $[-M,M] \subset \mathbb{R} $, EWA with the gradient trick can compete with the convex combination $\mathbf{q} \in \mathcal{X}$ of $N$ experts as shown by the regret bound \citep{devaine_forecasting_2013}:

\begin{equation}
    \mathcal{R}_T^{\text{EWA}^{grad}} (\mathbf{q}) \leqslant \frac{\ln(N)}{\eta} + \eta \frac{M^2}{2} T,
\end{equation}
with $T$ the number of iterations and $\eta >0$ the learning rate. To minimize this regret bound, the theoretical optimal choice for $\eta$ is $\eta^*= M^{-1} \sqrt{(2\ln(N))/T}$. With $\eta^*$ and if $\nabla \ell_t(\mathbf{w}_t) \cdot \mathbf{q} $ is bounded in $[-M,M] \subset \mathbb{R} $, the regret of EWA with the gradient trick against a convex combination $\mathbf{q}\in\mathcal{X}$ of experts $\mathcal{R}_T(\mathbf{q})$, becomes after $T$ iterations \citep{devaine_forecasting_2013}:

\begin{equation}
\label{eq:EWA_grad_optimal_regret}
    \mathcal{R}_T^{\text{EWA}^{grad}} (\mathbf{q}) \leqslant M \sqrt{2T \ln(N)}=o(T).
\end{equation}

Unfortunately, in both cases (with and without the gradient trick) $\eta^*$ depends on the number of iterations $T$ and the bound $M$, which are not always known in practice. Hence, $\eta^*$ is often not known in advance.

To remedy this, as suggested in \citet{stoltz_agregation_2010} and \citet{devaine_forecasting_2013}, one could run several EWA algorithms with different learning rates and choose at iteration $t$ the prediction of the EWA algorithm which has the smallest cumulative loss $\sum_{s=1}^{t-1} \ell_s(\mathbf{w}_s)$ at $t-1$. However, no theoretical guarantees are available for this method.

We did not pursue this method and obtained better results by using the adaptive learning rate of \citet{cesa-bianchi_improved_2007}:

\begin{equation}
    \eta_t= \min \left( {B_{t-1}}^{-1}, C\sqrt{\ln(N)/V_{t-1}} \right),
\end{equation}

where $t$ is the iteration, $N$ the number of experts, $C=\sqrt{2(\sqrt{2}-1/(e-2))}$, $E_t= \max \left( |\ell_{s}(\mathbf{\delta}_i)-\ell_{s}(\mathbf{\delta}_j)|, 1\leqslant s \leqslant t, 1 \leqslant i,j \leqslant N \right)$, $\ell_{t}(\mathbf{\delta}_i)$ is the loss of expert $i$, $B_t= 2^k$ with $k\in\mathbb{Z}$ the smallest integer such that $E_t \leqslant 2^k$, and $V_t=\sum_{s=1}^t \mathrm{Var}_s$ with $\mathrm{Var}_s=\sum_{i=1}^N w_{i,t}{\ell_{s}(\mathbf{\delta}_i)}^2 - \left( \sum_{i=1}^N w_{i,t}\ell_{s}(\mathbf{\delta}_i)\right)^2$ the variance of the loss under distribution $\mathbf{w}_t$. We used this learning rate in all of the following for EWA. This adaptive learning rate for EWA only comes with a cost of a factor smaller than 2 in the regret bound \citep{stoltz_agregation_2010}. For more details about the regret bound and the proof see \citet{cesa-bianchi_improved_2007}.

\subsection{Second order aggregations}
\label{second_order_part}
Aggregation strategies with second order regret bounds were introduced by \citet{cesa-bianchi_improved_2007}. They were improved by \citet{gaillard_second-order_2014} (MLprod and MLpol) with multiple learning rates. And later, \citet{koolen2015second} introduced Squint, followed by \citet{wintenberger_optimal_2017} with Bernstein Online Aggregation (BOA) and \citet{adjakossa_kalman_2023} with KAO.

The second order refinement consists in using the so called excess loss \citep{gaillard_second-order_2014} $\ell_{i,t}^{exc}= \ell_{t}(\mathbf{\delta}_i) - \sum_{i=1}^N w_{i,t} \ell_{t}(\mathbf{\delta}_i)$ instead of the loss $\ell_{t}(\mathbf{\delta}_i)$. So the excess loss $\ell_{i,t}^{exc}$ is the difference between the loss of expert $i$ and the weighted - by $\mathbf{w}_t$ - mean of the experts. It can be seen as an instantaneous regret of expert $i$ against the weighted mean of the experts.

These strategies with second order regret bounds allow one to take into account the difficulty of learning from the data and to perform better with "easy" data. For example, when the data has a small variance or when the experts have small losses \citep{koolen2015second,gaillard_second-order_2014}.

\citet{wintenberger_stochastic_2024} also showed that the second order refinement can provide nearly optimal regret bounds in a stochastic environment (where observations and experts are no longer black boxes but random variables).

For this reason, these strategies with second order regret bounds seem to be more adapted to weather forecasting. Indeed, the atmosphere will behave more like a stochastic process rather than an adversary who wants to make the loss of the aggregation as high as possible.

\begin{algorithm}
\caption{BOA aggregation \citep{wintenberger_stochastic_2024}}
\label{BOA_aggregation}
\begin{algorithmic}
    \State Initialization: $\mathbf{w_{1}} \in \mathbb{R}^N $ the uniform weight vector and  $\mathbf{\eta}_{0}=(0,\ldots,0)\in \mathbb{R}^N$ and $\tilde{L}_{i,0}=0$, $i=1,\ldots,N$.\\
    \For{$t \geqslant 1$}

    \begin{itemize}
        \item $\widehat{y}_t=\sum_{i=1}^N w_{i,t} x_{i,t}$
        \item $\ell_{i,t}^{exc}= \ell_{t}(\mathbf{\delta}_i) - \sum_{i=1}^N w_{i,t} \ell_{t}(\mathbf{\delta}_i)$ 
        \item $V_{i,t}=V_{i,t-1} + 2.2 {(\ell_{i,t}^{exc}})^2$, $i=1,\ldots,N$
        \item if $V_{i,t} \neq 0$, $\eta_{i,t}= \sqrt{1/V_{i,t}} $, otherwise $\eta_{i,t}=0$, $i=1,\ldots,N$
        \item surrogate loss: $\tilde{\ell}_{i,t}= \ell_{i,t}^{exc} + \eta_{i,t}(\ell_{i,t}^{exc})^2$, $i=1,\ldots,N$
        \item $\tilde{L}_{i,t}= \tilde{L}_{i,t-1} + \tilde{\ell}_{i,t}$, $i=1,\ldots,N$
        \item $w_{i,t+1}= \frac{w_{i,1} \eta_{i,t} \exp(-\eta_{i,t}\tilde{L}_{i,t})}{\sum_{j=1}^{N}w_{j,1} \eta_{j,t} \exp(-\eta_{j,t} \tilde{L}_{j,t})}$, $i=1,\ldots,N$
    \end{itemize}
    \EndFor
\end{algorithmic}
\end{algorithm}

The BOA algorithm in its adaptive version is presented in Algorithm \ref{BOA_aggregation}. Here, adaptive refers to the learning rates $\eta_{i,t}$ which depend on the expert $i=1,\ldots,N$ and on the iteration $t$. Moreover, computing them does not require prior knowledge of the loss bounds or $T$. \citet{wintenberger_stochastic_2024} showed that for BOA, the regret bound against expert $i=1,\ldots,N$ is\footnote{Let $f,g:\mathbb{R}\to\mathbb{R}$ be two functions with $g\ge0$. We say that $f$ is a big $O$ of $g$ and write $f(t)=O(g(t))$ if there exists a real number $C>0$ such that $|f(t)|\le Cg(t)$ for all $t\in\mathbb{R}$. We treat $ln(ln(T))$ as a constant.}:

\begin{equation}
\label{eq:regret_BOA_bound}
    \mathcal{R}_T^{\text{BOA}} (\mathbf{\delta}_i) \leqslant O \left( \sqrt{2.2 \sum_{t=1}^T (\ell_{i,t}^{exc})^2 } \right)=o(T),
\end{equation}
where $\ell_{i,t}^{exc}= \ell_{t}(\mathbf{\delta}_i) - \sum_{i=1}^N w_{i,t} \ell_{t}(\mathbf{\delta}_i)$ is the excess loss, $T$ the number of iterations, $\widehat{y}_t$ the EA's prediction, $y_t$ the observation and $x_{i,t}$ the prediction of expert $i$. This is a simplified version of the bound. For more details, cf. Theorem 4.1 of \citet{wintenberger_stochastic_2024}.

According to remark 6 of \citet{gaillard_second-order_2014}, BOA with the gradient trick has a regret against the fixed convex combination $\mathbf{q}=(q_1,\ldots,q_N)\in \mathcal{X}$ of experts:

\begin{equation}
\label{eq:regret_BOA_grad_bound}
    \mathcal{R}_T^{\text{BOA}^{grad}} (\mathbf{q}) \leqslant O \left( \sum_{i=1}^N q_i \sqrt{2.2 \sum_{t=1}^T \left(\frac{ \partial \ell(\widehat{y}_t,y_t)}{\partial \widehat{y}_t}(x_{i,t}-\widehat{y}_t)\right)^2 } \right)=o(T),
\end{equation}
where $T$ is the number of iterations, $\widehat{y}_t$ the EA's prediction, $y_t$ the observation, $x_{i,t}$ the prediction of expert $i$ and $\ell$ a convex differentiable loss function in its first argument with a bounded gradient.


\begin{algorithm}
\caption{MLprod aggregation \citep{gaillard_second-order_2014}}
\label{MLprod_aggregation}
\begin{algorithmic}
    \State Parameter: a rule to choose the learning rates $\eta_{i,t}$, for $i=1,\ldots,N$ and $t=0,\ldots,T$.
    \State Initialization: initial weight vector $\mathbf{w}_{1}\in \mathbb{R}^N$ the uniform weight vector.\\
    \For{$t \geqslant 1$}

    \begin{itemize}
        \item $\widehat{y}_t=\sum_{i=1}^N w_{i,t} x_{i,t}$
        \item observe the losses $\ell_t(\mathbf{w}_t)$ and $\ell_{t}(\mathbf{\delta}_i)$
        \item $\ell_{i,t}^{exc}= \ell_{t}(\mathbf{\delta}_i) - \sum_{i=1}^N w_{i,t} \ell_{t}(\mathbf{\delta}_i)$ 
        \item according to the rule, choose $\eta_{i,t}>0$, $i=1,\ldots,N$
        \item $ w_{i,t+1}= \frac{ \eta_{i,t} \left( w_{i,t} (1 - \eta_{i,t-1} \ell_{i,t}^{exc}) \right)^ \frac{\eta_{i,t}}{\eta_{i,t-1}} }{ \sum_{j=1}^N\eta_{j,t} \left( w_{j,t} (1 - \eta_{j,t-1} \ell_{j,t}^{exc}) \right)^ \frac{\eta_{j,t}}{\eta_{j,t-1}} }$, $i=1,\ldots,N$
    \end{itemize}
    
    \EndFor
\end{algorithmic}
\end{algorithm}

The MLprod algorithm is presented in Algorithm \ref{MLprod_aggregation}. It is also an adaptive algorithm such as BOA, presented in Algorithm \ref{BOA_aggregation}. \citet{gaillard_second-order_2014} showed that if the losses $\ell_t(\mathbf{\delta}_i)$ are bounded in $[0,1]$ and the learning rates are such that $\eta_{i,t}=\min \{0.5 , \sqrt{\ln(N)/(1+V_{i,t})} \}$, with $V_{i,t}=\sum_{s=1}^t ({\ell_{i,s}^{exc}})^2$, where $\ell_{i,t}^{exc}= \ell_{t}(\mathbf{\delta}_i) - \sum_{i=1}^N w_{i,t} \ell_{t}(\mathbf{\delta}_i)$ is the excess loss, and $\ell_t(\mathbf{\delta}_i)$ the loss of expert $i=1,\ldots,N$, $t=1,\ldots,T$, then MLprod has a regret bound against expert $i$:
\begin{equation}
\label{eq:regret_MLprod_bound}
    \mathcal{R}_T^{\text{MLprod}}(\mathbf{\delta}_i) \leqslant O \left( \frac{1}{\sqrt{\ln(N)}}\sqrt{1+ \sum_{t=1}^T (\ell_{i,t}^{exc})^2} + \ln(N) \right)=o(T).
\end{equation}

This is a simplified bound, for more details cf. Corollary 4 of \citet{gaillard_second-order_2014}.

According to remark 6 of \citet{gaillard_second-order_2014}, the regret of MLprod with the gradient trick against the fixed convex combination $\mathbf{q}=(q_1,\ldots,q_N)\in \mathcal{X}$ of experts is:

\begin{equation}
\label{eq:regret_MLprod_grad_bound}
    \mathcal{R}_T^{\text{MLprod}^{grad}}(\mathbf{q}) \leqslant O \left( \sum_{i=1}^N q_i \frac{1}{\sqrt{\ln(N)}}\sqrt{1+ \sum_{t=1}^T \left(\frac{ \partial \ell(\widehat{y}_t,y_t)}{\partial \widehat{y}_t}(\widehat{y}_t-x_{i,t})\right)^2} + \ln(N) \right)=o(T),
\end{equation}

where $\ell$ is the convex differentiable loss function in its first argument with a bounded gradient, $\widehat{y}_t$ the EA's prediction, $y_t$ the observation, $x_{i,t}$ the prediction of expert $i$, $N$ the number of experts and $T$ the number of iterations.

\begin{algorithm}
\caption{MLpol aggregation \citep{gaillard_second-order_2014}}
\label{MLpol_aggregation}
\begin{algorithmic}
    \State Parameter: a rule for choosing the learning rates $\eta_{i,t}$, for $i=1,\ldots,N$ and $t=1,\ldots,T$.
    \State Initialization: $\mathbf{\mathcal{R}}_0^{exc}=(0,\ldots,0) \in \mathbb{R}^N $,  $\mathbf{w}_1 \in \mathbb{R}^N$ the uniform weight vector.\\
    \For{$t \geqslant 1$}

    \begin{itemize}
        \item $\widehat{y}_t=\sum_{i=1}^N w_{i,t} x_{i,t}$
        \item observe the losses $\ell_{t}(\mathbf{w}_t)$ and $\ell_{t}(\mathbf{\delta}_i)$, $i=1,\ldots,N$
        \item $\mathcal{R}_{i,t}^{exc}=\mathcal{R}_{i,t-1}^{exc}-\ell_{i,t}^{exc}$, $i=1,\ldots,N$
        \item according to the rule, choose $\eta_{i,t}$, $i=1,\ldots,N$
        \item $w_{i,t+1}= \frac{\eta_{i,t}(\mathcal{R}_{i,t}^{exc})_{+}}{\sum_{j=1}^N \eta_{j,t}(\mathcal{R}_{j,t}^{exc})_{+}}$, $i=1,\ldots,N$, where $x_{+}$ denotes the vector of the non negative parts of the components of $x$
    \end{itemize}
    
    \EndFor
\end{algorithmic}
\end{algorithm}

The adaptive MLpol algorithm is presented in Algorithm \ref{MLpol_aggregation}. \citet{gaillard_second-order_2014} showed that if the loss $\ell_t$ is bounded between 0 and 1 and if $\eta_{i,t}=(1+V_{i,t})^{-1}$, with $V_{i,t}=\sum_{s=1}^t ({\ell_{i,s}^{exc}})^2$ and $\ell_{i,t}^{exc}= \ell_{t}(\mathbf{\delta}_i) - \sum_{i=1}^N w_{i,t} \ell_{t}(\mathbf{\delta}_i)$ the excess loss, $i=1,\ldots,N$, $t=1,\ldots,T$, then MLpol has a regret bound against expert $i$:

\begin{equation}
\label{eq:regret_MLpol_bound}
    \mathcal{R}_T^{\text{MLpol}}(\mathbf{\delta}_i) \leqslant O \left( \sqrt{ N(1+\ln(T)) \left(1+\sum_{t=1}^T (\ell_{i,t}^{exc})^2\right) } \right)=o(T),
\end{equation}

and with the gradient trick, a regret bound against the convex combination $\mathbf{q}=(q_1,\ldots,q_N) \in \mathcal{X}$ of experts:
\begin{equation}
\label{eq:regret_MLpol_grad_bound}
    \mathcal{R}_T^{\text{MLpol}^{grad}}(\mathbf{q}) \leqslant O \left( \sum_{i=1}^N q_i \sqrt{ N(1+\ln(T)) \left(1+\sum_{t=1}^T \left(\frac{ \partial \ell(\widehat{y}_t,y_t)}{\partial \widehat{y}_t}(\widehat{y}_t-x_{i,t})\right)^2\right) } \right)=o(T),
\end{equation}

where $\ell$ is the convex differentiable loss function in its first argument with a bounded gradient.

To back up that aggregations with second order bounds are adapted for non adversarial settings especially when the losses are "small", it can be observed that if $\ell$ is the square loss we have:
\begin{equation}
    \left( \frac{ \partial \ell(\widehat{y}_t,y_t)}{\partial \widehat{y}_t} \right)^2=(2(\widehat{y}_t - y_t))^2=4\ell_t(\mathbf{w}_t).
\end{equation}
So by replacing ${\frac{ \partial \ell(\widehat{y}_t,y_t)}{\partial \widehat{y}_t}}^2$ by $4\ell_t(\mathbf{w}_t)$ in the second order regret bounds, one can see that when the aggregations error is small the regret is small too.

\section{Data and evaluation}
\label{section_data_evaluation}
\subsection{Data and experts}
\label{subsection_data_experts}

\begin{figure}[h]
 \centerline{\includegraphics[width=27pc]{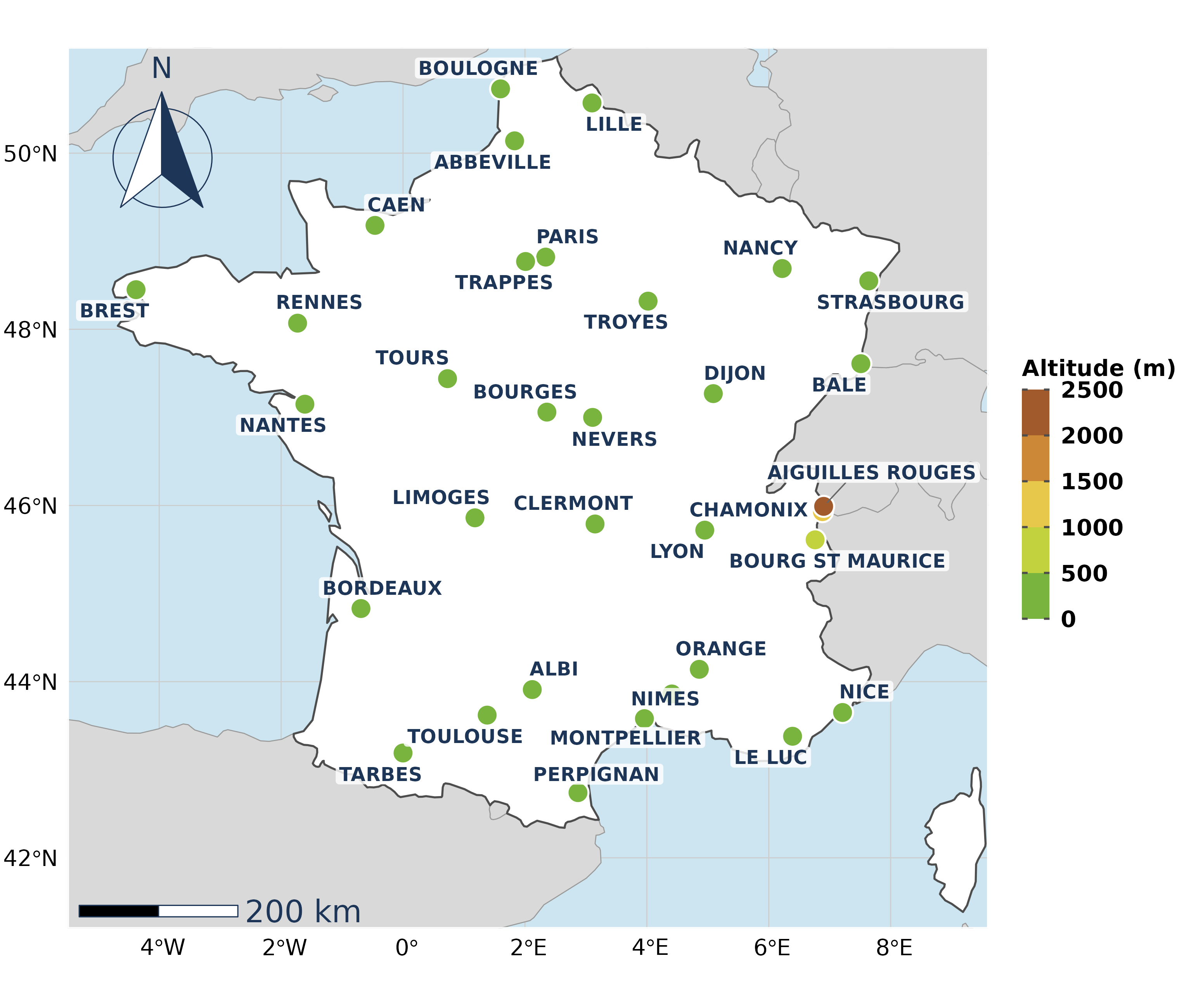}}
  \caption{Map of meteorological station locations colored by altitude across France (in white). All stations are situated at low elevations ($\leqslant$500 m) except for the stations in the French Alps: Aiguilles Rouges (brown dot) and Chamonix (yellow dot) which are in the same valley and Bourg St Maurice (light green dot). Latitude on the ordinate, longitude on the abscissa. The sea is in blue and the neighboring states are in gray.}\label{map_stations}
\end{figure}

The data used in this study are available through the link provided in the data availability statement before the Appendix. 
They cover the period from 2020-03-30 to 2023-09-03. We used the temperature observations 2 meters above the ground in 33 stations located on the map in Figure \ref{map_stations}. The chosen stations are professional weather stations, providing good spatial coverage as well as a good representation of the various climate conditions in France. All stations have an uncertainty of 0.5°C, except Aiguilles rouges, for which the uncertainty is 1°C. 

We also have the predictions for the base time 00 a.m UTC of several raw NWP models: forecasts of Application of Research to Operations at Mesoscale (AROME), Integrated Forecasting System (IFS) of the European Centre for Medium-Range Weather Forecasts (ECWMF) and Action de Recherche Petite Echelle Grande Echelle (ARPEGE). They are respectively denoted by raw.aro, raw.arp and raw.ifs. 

ARPEGE and IFS are hydrostatic global NWP models. ARPEGE has a varying horizontal resolution from 5km over France to 24 km at the antipodes. The horizontal resolution of IFS is 9km. AROME is a non-hydrostatic limited area NWP model over Europe and has a horizontal resolution of 1.2 km.

In addition, we use post-processed forecasts of AROME, ARPEGE and IFS. They are denoted by ppm.aro, ppm.arp and ppm.ifs. The post-processing of AROME, ARPEGE and IFS is done at station points using random forests \citep{breiman_random_2001} which are simple non parametric machine learning methods.

 Finally, we used the predictions of the 10\%, 30\%, 50\%, 70\% and 90\% quantiles of the post processed Prévision d’Ensemble ARPEGE (PEARP) \citep{taillardat_research_2020} denoted respectively as Q10, Q30, Q50, Q70 and Q90. The post processed PEARP quantiles are defined station wise, thanks to quantile random forests \citep{taillardat_research_2020}.

The temperature predictions of the raw NWP models are bilinearly interpolated on the stations. For the PPMs and the post processed quantiles, a bilinear interpolation is done for different variables on the stations without taking into account the altitude. Then, at each station, post-processing is done with random forests, which account for errors introduced by inaccurate model elevations. 

The lead times range from 6 to 48 hours with a 3 hours step for all the models plus the lead times 57, 72 and 84 hours for the raw and PPMs ARPEGE, IFS and the PEARP quantiles.

\subsection{Benchmarks}
\label{benchmarks}
The first benchmarks are the different experts. We want to know if the EAs beat the "best expert" oracle. We call "best expert", expert $i^* \in \{1,\ldots,N\}$ such that $i^*=\arg \min_{i=1,\ldots,N} \sum_{t=1}^T \ell_{t}(\mathbf{\delta}_i)$ for a specific lead time and station. It is the best possible prediction on average that one can achieve by choosing only one expert from start to end among experts $i=1,\ldots,N$. Hence, trivially, if the EA is better than the best expert, then it outperforms on average all aggregated raw and post-processed NWP models.

The second benchmark is the "best convex combination" of experts oracle. Analogously to the best expert, this is the strategy which predicts like $\mathbf{q}^*\in \mathcal{X}$ the best convex combination of experts in hindsight for a specific station and lead time: $\mathbf{q}^*=\arg \min_{\mathbf{q}\in \mathcal{X}} \sum_{t=1}^T \ell_{t}(\mathbf{q})$. It is the best possible prediction on average that one can achieve by choosing a fixed - from start to end - convex combination among experts $i=1,\ldots,N$.

The third benchmark is the operational NBM model of \citet{craven_national_2020}. It is a two step strategy. At $t\geqslant1$, the first step is to reduce the bias of all experts $i=1,\ldots,N$:

\begin{equation}
    x_{i,t}^b=x_{i,t}-B_{t}
\end{equation}
with

\begin{equation}
B_t=(1-\alpha)B_{t-1}+\alpha(x_{i,t-1}-y_{t-1})
\end{equation}

where $\alpha$ is the decaying weight factor. It plays a similar role to that of the learning rates in EA. In NBM v3.2, $\alpha=0.025$ from September to May and 0.05 otherwise.
The second step is to make a weighted average, and to predict $\widehat{y}_t=\sum_{i=1}^N w_{i,t}x_{i,t}^b$ with:

\begin{equation}
    w_{i,t}=\frac{{L_{i,t}}^{-1}}{\sum_{j=1}^N {L_{j,t}}^{-1}}
\end{equation}

and

\begin{equation}
    L_{i,t}=(1-\alpha)L_{i,t-1}+\alpha \ell_{t-1}(\mathbf{\delta}_i)
\end{equation}

for $t\geqslant1$, $\alpha$ the decaying weight factor and $\ell_t(\mathbf{\delta}_i)$ the loss of expert $i$. In \citet{craven_national_2020}, the loss is the absolute loss, but for fairness, we replaced it by the square loss since we are using the RMSE to evaluate the aggregations.

The major difference between NBM and EA, is that NBM tries to improve the experts by removing the bias before averaging them. Whereas EA only tries to choose the best expert or the best convex combination of experts, without changing the experts. Furthermore, the decaying weight factor $\alpha$ is not adaptive, whereas the adaptive learning rates that we used for EA depend on the iteration and on the experts.

\subsection{Evaluation methods}
To compare the aggregation strategies and benchmarks across all stations and lead times, we used several scoring functions since they have different purposes \citep{gneiting_making_2011}.

Our main scoring function is the Root Mean Squared Error (RMSE). We recall that the RMSE for $k$ predictions $x_1,\ldots,x_k$ and observations $y_1,\ldots,y_k$ is $\sqrt{k^{-1} \sum_{s=1}^k (x_s - y_s)^2 }$. The RMSE evaluates the ability of the forecast to predict the mean of a distribution and it is sensitive to outliers since it penalizes more large errors. This last point is quite convenient for us since we consider large temperature errors to be worse than small errors.

That is why we also looked at the 95\% quantile of the absolute error denoted $Q_{95}(|e|)$ because minimizing the $Q_{95}(|e|)$ reduces large errors. Indeed, the $Q_{95}(|e|)$ is only sensitive to very large errors and can be considered as an "extremes metric". Hence the RMSE and the $Q_{95}(|e|)$ are complementary. 

And finally, we also computed the bias which indicates if the error is centered in 0 or if there is a systematic error. For $k$ predictions $x_1,\ldots,x_k$ and observations $y_1,\ldots,y_k$, the bias is $k^{-1}\sum_{s=1}^k (x_s - y_s)$.

In order to assess if our final results of Section \ref{section_results} are significant we did some statistical tests. For the RMSE differences, we did a Diebold–Mariano (DM) test \citep{diebold_comparing_1995} as in \citet{schulz_machine_2022} for every couple station-lead time. We also did quantile tests \citep{johnson_two-sample_1987} for the $Q_{95}(|e|)$ differences. As advised by \citet{wilks_stippling_2016}, we did a Benjamini-Hochberg (BH) procedure \citep{benjamini_controlling_1995} for both tests like in \citet{rasp_neural_2018,schulz_machine_2022,pic_mathematical_2022} to handle multiple testing with spatial and temporal dependencies for a 5\% level. The results of these tests can be found in the Appendix (Table \ref{test_rmse} and \ref{test_q95}).

\section{Results}
\label{section_results}

In order to run the EAs we used the Opera \citep{gaillard_opera_2016} package (version 1.1.1) with R (version 3.6.1). The code can be found via the link in the data availability statement at the end of the paper. We always run one EA for each station-lead time pair.

\subsection{Influence of the set of experts}
\label{section_influence_number_experts}

We now study how BOA reacts to different sets of experts. The following equation of \citet{bourel_boosting_2024} decomposes the aggregation's squared error in the mean of the expert's error (a) and the variance/diversity of the experts (b):
\begin{equation}
    (\widehat{y_{t}}-y_{t})^{2}=\underbrace{ \sum_{j=1}^{N}w_{j,t}(x_{j,t}-y_{t})^{2} }_{(a)} - \underbrace{ \sum_{j=1}^{N}w_{j,t}(x_{j,t}-\widehat{y}_{t})^{2} }_{(b)},
\label{bias_variance_tradeoff_square_loss}
\end{equation}
with $\widehat{y}_t$ the EA's prediction, $y_t$ the observation and $w_{j,t}$ and $x_{j,t}$ respectively the weight and the prediction of expert $j$ at $t$.

This equation reflects the fact that the square loss of the aggregation increases with the loss of the experts and decreases with the diversity of the experts. Furthermore, \citet{gaillard_forecasting_2015} wrote that "[...] aggregating rules are quite robust to adding experts, and [...] combining forecasts does not suffer much from over fitting". So, adding uncorrelated models with similar accuracy is likely to improve the aggregation.

\begin{figure*}[h]
 \centerline{\includegraphics[width=39pc]{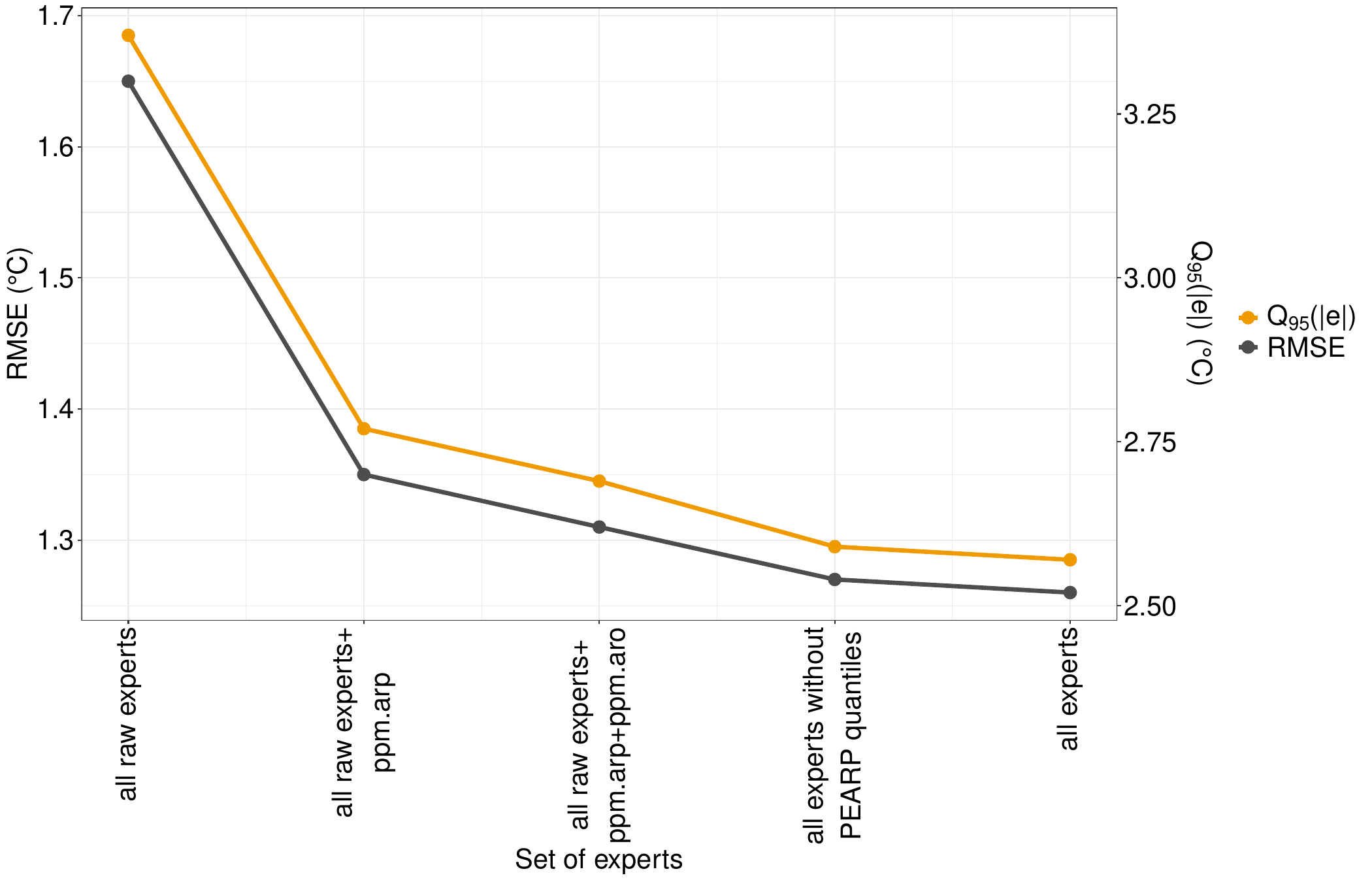}}
  \caption{Scores of the BOA \citep{wintenberger_stochastic_2024} aggregation (without the gradient trick) depending on the available experts. RMSE (black) on the left ordinate and 95 \% quantile of the absolute error $Q_{95}(|e|)$ (orange) on the right ordinate. On the abscissa, the sets of experts. From the left to the right, there are more and more experts.}\label{influence_number_experts}
\end{figure*}

In Figure \ref{influence_number_experts}, we can see the RMSE and $Q_{95}(|e|)$ of BOA depending on the available experts. The scores become better when the number of experts grows. Especially adding the PPMs improves the scores while the improvement due to the PEARP quantiles is quite low. 

We only show the results for BOA, but we observed the same behavior for the other EAs, except when we added the PEARP quantiles. As we added them, BOA and EWA improved whereas Mlpol kept the same scores and MLprod was slightly worse. This may be due to a negligible difference and to the fact that EWA and BOA can adapt to the range of losses.

Thus, it is not obvious that our EAs are able to deal with these biased experts, which are bad on average and could add some noise to the aggregation.

\begin{figure*}[h]
 \centerline{\includegraphics[width=39pc]{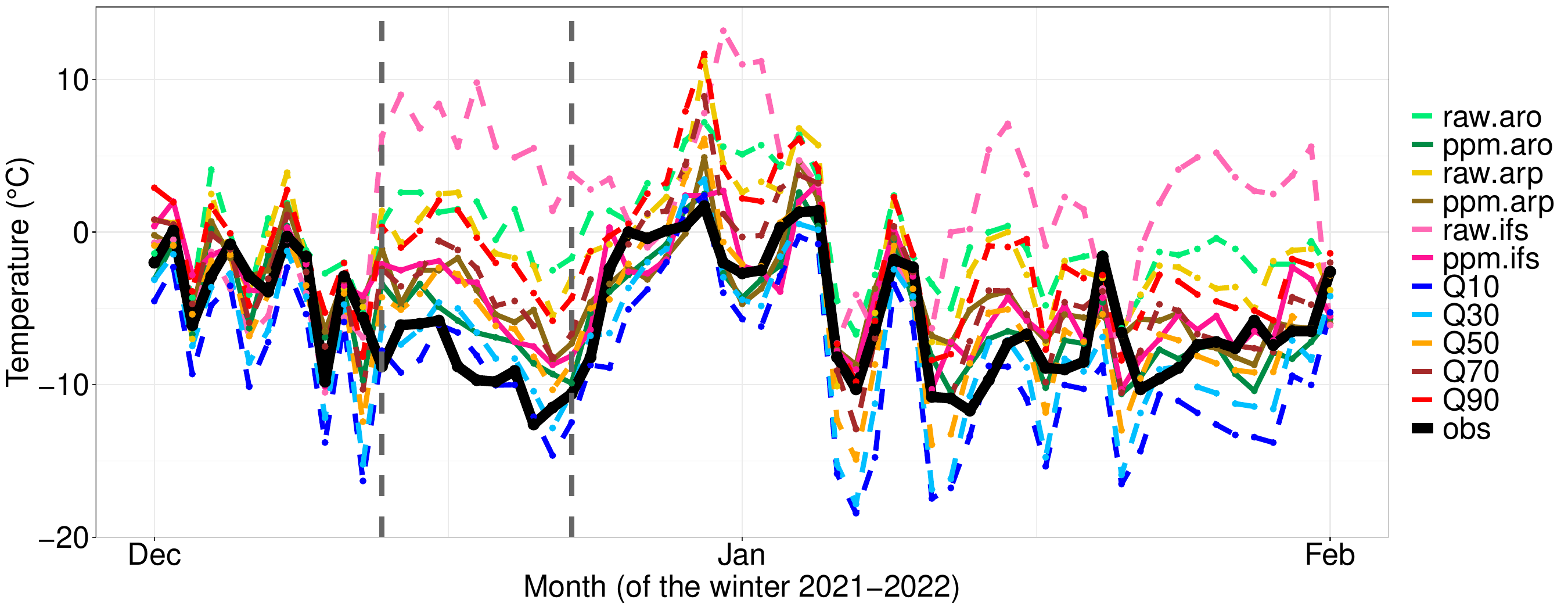}}
  \caption{Observed temperatures (obs) in black and predicted temperatures for the Chamonix station (French Alps) for the lead time 48 hours of the midnight run, during winter 2021-2022. One can see the predictions of the "raw" and Post Processed Models (ppm) AROME (aro), IFS (ifs), ARPEGE (arp) and of the post processed quantiles of the PEARP modele (Q10, Q30, Q50, Q70, Q90). The dashed grey vertical lines delimit the mid-December period where all the experts make too warm predictions except the Q10 and Q30 experts.}\label{chamonix_obs_all_experts}
\end{figure*}

However, Figure \ref{chamonix_obs_all_experts} highlights that without the PEARP quantiles, the observations would often be outside of the experts’ prediction set. For instance, during a ten-day period in mid-December 2022 (indicated by the dashed grey lines), the Q10 and Q30 experts were the only experts to remain close to the observations in Chamonix.

As commented in Subsection \ref{section_theoretical_framework}.\ref{subsection:expert_agregation}, we use EA strategies which are convex and hence cannot predict an observation that falls outside the range of the experts' predictions. Therefore, we added the 10\% ,30\%, 70\% and 90\% post-processed PEARP quantiles to the set of experts so that the observation almost always lies in the range of the experts' predictions.

The inclusion of all experts is also justified by Equation \eqref{bias_variance_tradeoff_square_loss} as explained before, and because the regret bounds of the EAs are weakly dependent on the number of experts.

Hence, the theory and our experiments suggest that we should use a large number of experts in order to make the best predictions.

\subsection{Can we beat the best expert in hindsight?}
\label{subsection:can_we_beat_the_best_expert_in_hindsight}

As seen previously in Section \ref{section_theoretical_framework}, the theory shows that our EAs - without the gradient trick - should be able to be asymptotically as good as (or even better than) the best expert in hindsight since their regret satisfies the condition $\mathcal{R}_T\le o(T)$. The scores of the EAs, without the gradient trick, with all the raw models, PPMs and PEARP quantiles are reported in Table \ref{rmse_ewa_boa_MLpol_without_grad_with_pearp}.

The bias is small for every aggregation strategy except for the uniform aggregation which has a slightly larger negative bias.

With regard to all scores, all EAs outperform the uniform aggregation but all are worse than the best convex combination of experts. Only BOA, EWA and MLprod have a smaller RMSE and $Q_{95}(|e|)$ than the best expert, hence outperform all raw models and PPMs.

It is quite surprising that EWA performs better than MLprod and MLpol, despite their second order regret bound. This shows that having the regret bounds of the aggregation strategies - which can sometimes be difficult to compare - is not enough to conclude which of the aggregation strategies should be the best. So in practice, one has to test the different aggregation strategies in order to estimate which strategy performs best.

\begin{table}[h]
\centering
\begin{tabular}{|c|c|c|c|}
\hline
aggregation & BIAS (°C) & RMSE (°C) & $Q_{95}(|e|)$\\
\hline
EWA & -0.01 & 1.27 & 2.58 \\
BOA & -0.01 & 1.26 & 2.57 \\
MLpol & -0.01 & 1.38 & 2.8 \\
MLprod & -0.01 & 1.31 & 2.66 \\
uniform aggregation & -0.03 & 1.41 & 2.90 \\
best Expert & -0.01 & 1.37 & 2.8 \\
best convex combination & -0.01 & 1.23 & 2.5 \\
\hline
\end{tabular}
\caption{Bias, RMSE and 95\% quantile of the absolute error $Q_{95}(|e|)$ for all the lead times and stations. The scores are for the uniform aggregation, the best expert oracle and the best convex combination of experts oracle, and the aggregations (without the gradient trick and with all the experts): EWA \citep{vovk_aggregating_1990}, BOA \citep{wintenberger_stochastic_2024}, MLprod and MLpol \citep{gaillard_second-order_2014}. The best expert is the strategy predicting like the best (on average) expert in hindsight for every couple station/lead time. The best convex combination is the strategy predicting like the best (on average) fixed convex combination of experts in hindsight for every couple station/lead time.}
\label{rmse_ewa_boa_MLpol_without_grad_with_pearp}
\end{table}

\begin{figure*}[h]
 \centerline{\includegraphics[width=39pc]{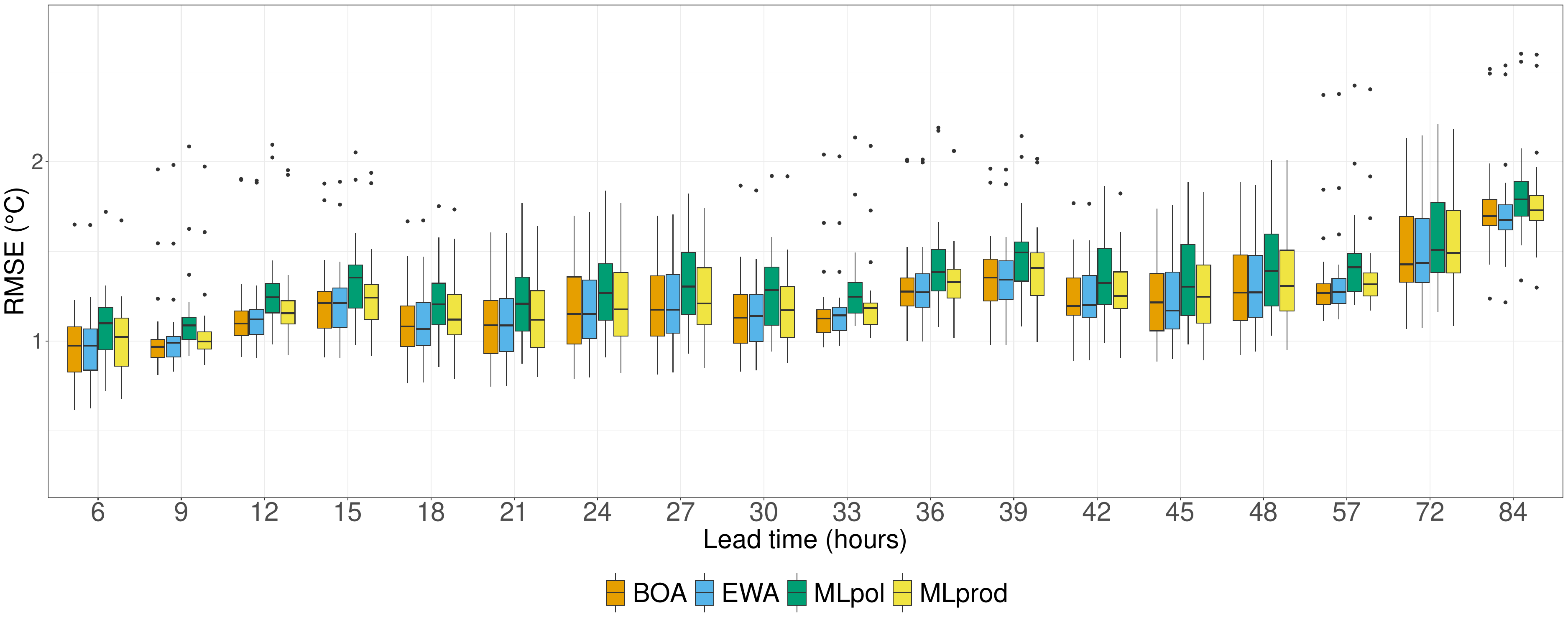}}
  \caption{Box plots of station-level RMSE across lead times for the following aggregation methods: BOA \citep{wintenberger_stochastic_2024}, EWA \citep{vovk_aggregating_1990}, MLpol and MLprod \citep{gaillard_second-order_2014} with all experts but without the gradient trick. Each box summarizes the spread across stations at a given lead time.}\label{boxplot_with_pearp_without_grad}
\end{figure*}

In Figure \ref{boxplot_with_pearp_without_grad} we can see the box plots of the RMSE of the stations for different aggregations (without the gradient trick and with all the experts) and lead times. It shows that BOA and EWA have similar distributions and are better than MLprod and MLpol.

We may also observe that the RMSE median becomes larger when the lead time increases but it is not completely monotonic since the box plots also show that there are more errors during the day than during the night. This is most likely due to the fact that the temperature is more predictable during the night because there is less convection and the atmosphere is more stable. We can also see that the RMSE variability depends on the the lead time.

Thus, our EAs without the gradient trick enable us to outperform the best expert oracle (hence all aggregated experts). But these EAs are outperformed by the best convex combination of experts oracle.

\subsection{Can we beat the best convex combination in hindsight?}

With the gradient trick, we are able to compete with the best convex combination of experts (see inequalities (\ref{eq:EWA_grad_optimal_regret}), (\ref{eq:regret_BOA_grad_bound}), (\ref{eq:regret_MLprod_grad_bound}) and (\ref{eq:regret_MLpol_grad_bound})). At first glance, this is better than competing with the best expert, because the best convex combination of experts is at least as good as the best expert. However, it is not certain that using the gradient trick will improve the results in practice. For example, the regret bound of EWA against convex combinations of experts is worse than against the experts themselves as equations \eqref{eq:EWA_optimal_regret} and \eqref{eq:EWA_grad_optimal_regret} show.

Therefore, we performed the same experiment as in subsection \ref{section_results}.\ref{subsection:can_we_beat_the_best_expert_in_hindsight}, but this time with the gradient trick and all experts. In Table \ref{rmse_ewa_boa_MLpol_with_grad_with_pearp} we can see the scores obtained for these experiments. The biases are unchanged but the RMSE and the $Q_{95}(|e|)$ of all EAs are improved thanks to the gradient trick. All EAs have better scores than the best expert oracle. MLpol is improved to the extent that it has (along with BOA) the best scores of all EAs.

We can also see that the aggregations are much closer (compared to the EAs without the gradient trick in Table \ref{rmse_ewa_boa_MLpol_without_grad_with_pearp}) to the RMSE and $Q_{95}(|e|)$ of the best convex combination in hindsight. In particular, MLpol and BOA have almost the same RMSE as this oracle. 
So, the gradient trick improves the EAs, and enables us to approach the scores of the best convex combination of experts. But it remains to know if these scores are statistically different or not.

The corresponding quantile tests with a BH procedure are presented in the Appendix Table \ref{test_q95}. They show that the differences between the $Q_{95}(|e|)$ of BOA, MLpol, Mlprod and EWA are not statistically significant with a 0.05 level. It also shows that all these EAs have a significantly better $Q_{95}(|e|)$ for at least 18\% of the station-lead time pairs compared to the uniform aggregation and for at least 37\% compared to the best expert oracle.

In Table \ref{test_rmse} of the Appendix, we can see the DM tests with a BH procedure. It reveals that with the gradient trick and all experts, all the EAs significantly outperform the RMSE of the best expert for all station-lead time pairs. Furthermore, BOA is significantly better for 0.2\% of the station-lead time pairs than EWA and MLprod. And finally, none of the EAs outperform the best convex combination of experts.

Hence, with the gradient trick, all EAs largely outperform all the experts including the raw and post-processed IFS for example, which makes it interesting to use EA operationally. But the best convex oracle remains significantly better than all EAs.

\begin{table}[h]
\centering
\begin{tabular}{|c|c|c|c|}
\hline
aggregation & BIAS (°C) & RMSE (°C) & $Q_{95}(|e|)$\\
\hline
RAW ARP & -0.33 & 2.27 & 4.4 \\
PPM ARP & -0.01 & 1.57 & 2.9 \\
RAW IFS & 0.03 & 1.95 & 4.0 \\
PPM IFS & -0.01 & 1.42 & 2.9 \\
NBM & -0.01 & 1.31 & 2.64\\
EWA & -0.01 & 1.25 & 2.54\\
BOA & -0.01 & 1.24 & 2.53\\
MLpol & -0.01 & 1.24 & 2.53\\
MLprod & -0.01 & 1.25 & 2.54\\
uniform aggregation & -0.03 & 1.41 & 2.9\\
best Expert & -0.01 & 1.37 & 2.8\\
best convex combination & -0.01 & 1.23 & 2.5\\
\hline
\end{tabular}
\caption{Bias, RMSE and 95\% quantile of the absolute error $Q_{95}(|e|)$ for all the lead times and stations. The scores are for the uniform aggregation, the best expert oracle and the best convex combination of experts oracle, the "raw" and post processed models (PPM) ARPEGE (ARP) and IFS and the aggregations (with the gradient trick and with all the experts): EWA \citep{vovk_aggregating_1990}, BOA \citep{wintenberger_stochastic_2024}, MLprod and MLpol \citep{gaillard_second-order_2014}. The best expert is the strategy predicting like the best (on average) expert in hindsight for every couple station/lead time. The best convex combination is the strategy predicting like the best (on average) fixed convex combination of experts in hindsight for every couple station/lead time.}
\label{rmse_ewa_boa_MLpol_with_grad_with_pearp}
\end{table}

\subsection{Can we beat NBM?}

It now remains to compare the EAs with NBM. As we can see in Table \ref{rmse_ewa_boa_MLpol_with_grad_with_pearp}, NBM has a larger $Q_{95}(|e|)$ than BOA, EWA, MLprod and MLpol. But the statistical tests of Table \ref{test_q95} in the Appendix show that these differences are not statistically significant.

Table \ref{rmse_ewa_boa_MLpol_with_grad_with_pearp} also shows that NBM has a larger RMSE than the EAs. But this time, Table \ref{test_rmse} of the Appendix shows that the differences are statistically significant for 93.3\% of the station-lead time pairs in favor of BOA and a little bit less in favor of MLpol, EWA and MLprod. 

\begin{figure*}[h]
 \centerline{\includegraphics[width=39pc]{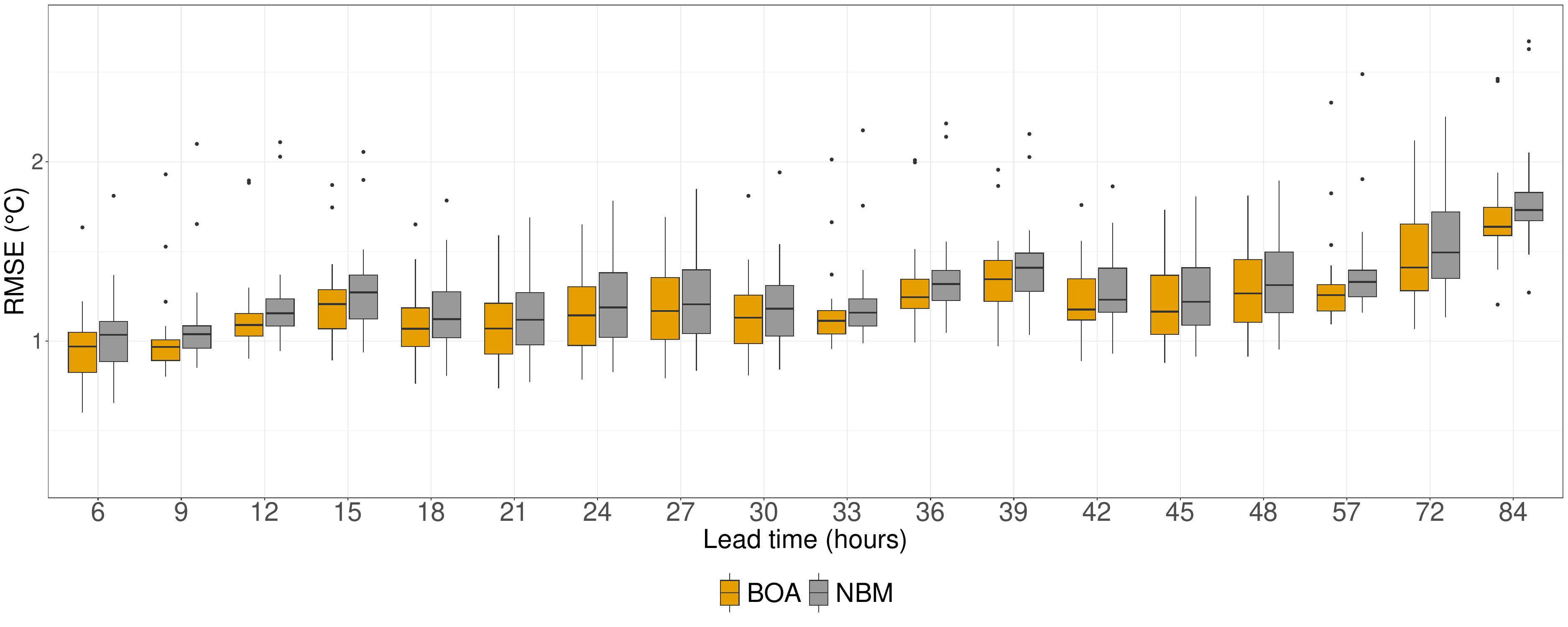}}
  \caption{Box plots of station-level RMSE across lead times for two aggregation methods: National Blend of Models (NBM) \citep{craven_national_2020} and Bernstein Online Aggregation (BOA) \citep{wintenberger_stochastic_2024} with the gradient trick and all experts. Each box summarizes the spread across stations at a given lead time.}\label{boxplot_with_pearp_with_grad}
\end{figure*}

In Figure \ref{boxplot_with_pearp_with_grad} we can see the box plots of the RMSE of the stations for different lead times, of NBM and BOA with the gradient trick and all the experts. For clarity, we did not plot the boxplots of EWA, MLpol and MLprod because they are very similar to BOA's boxplots when the gradient trick is used.

These boxplots show that the gradient trick has a huge impact on the aggregation. For every lead time the RMSE of BOA is improved compared to Figure \ref{boxplot_with_pearp_without_grad}. Again, we can see that the scores are better at night-time lead times.

It also shows that BOA has a better RMSE distribution than NBM for every lead time reinforcing that BOA is better than NBM. So NBM is outperformed by the EAs with the gradient trick, and it is statistically significant for up to 93.3\% of the station-lead time pairs in favor of BOA.

\subsection{Weight behavior}
\begin{figure*} [h]
 \centerline{\includegraphics[width=39pc]{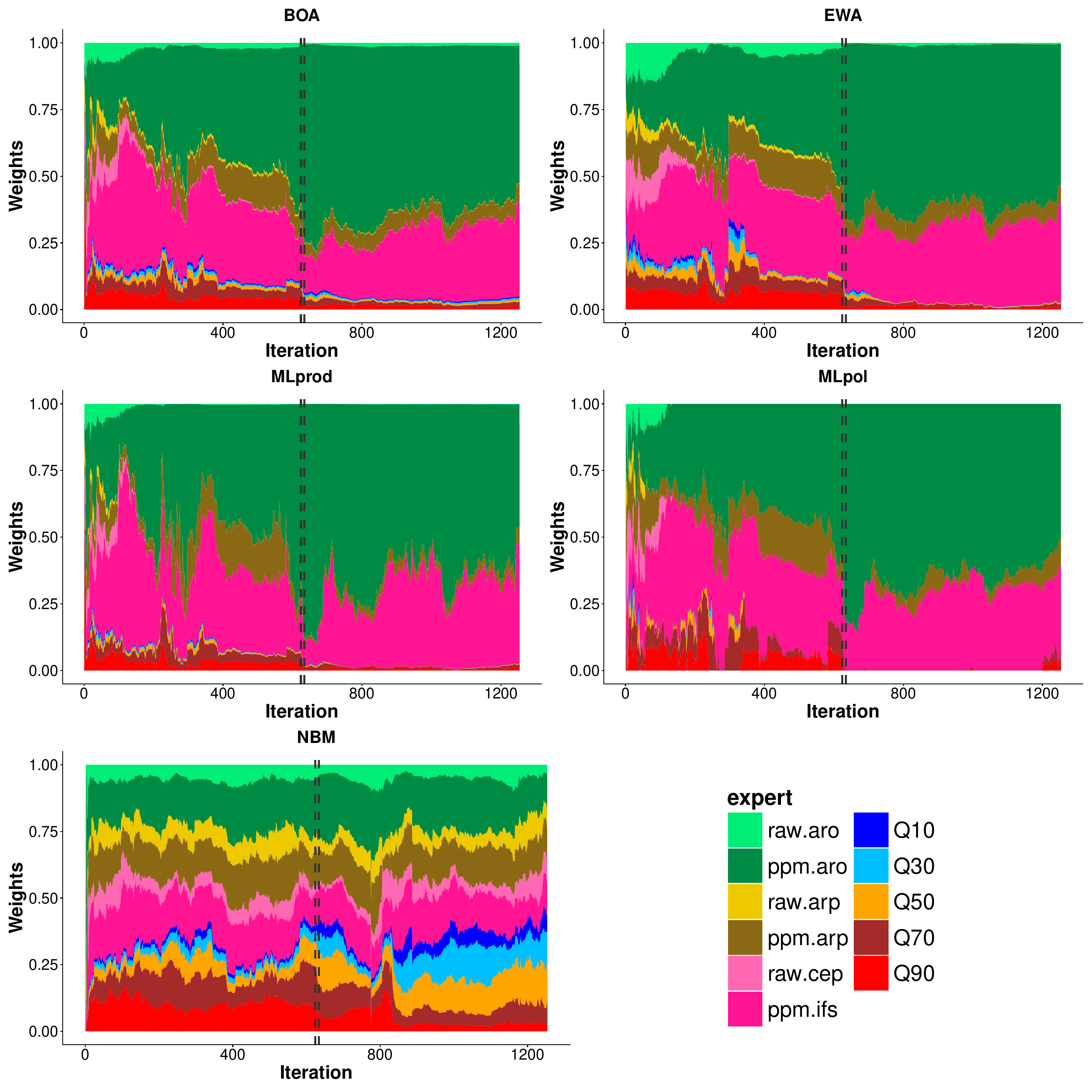}}
  \caption{Weights of the experts in Chamonix, from 2020-03-30 to 2023-09-03, for the lead time 48 hours, for NBM \citep{craven_national_2020} and different expert aggregations with the gradient trick (EWA \citep{vovk_aggregating_1990}, BOA \citep{wintenberger_stochastic_2024}, MLpol and MLprod \citep{gaillard_second-order_2014}). The experts are the "raw" and Post Processed Models (ppm) AROME (aro), IFS (ifs), ARPEGE (arp) and the post processed quantiles of the PEARP modele (Q10, Q30, Q50, Q70, Q90). The two dashed grey - almost overlapping - lines delimit the mid-December period where all experts make too warm predictions except the Q10 and Q30 experts. }\label{BOA_MLpol_MLprod_EWA_window1253_ech48_74056001_2020-03-30_2023-09-03_weights}
\end{figure*}

Now let us observe the weights of the experts given by the EAs and NBM. Figure \ref{BOA_MLpol_MLprod_EWA_window1253_ech48_74056001_2020-03-30_2023-09-03_weights} shows that for all the EAs, the weights of the biased experts are very quickly negligible and become a sort of noise at the end.

Figure \ref{BOA_MLpol_MLprod_EWA_window1253_ech48_74056001_2020-03-30_2023-09-03_weights} also shows that the different EAs have trouble using the biased experts even when only some of the biased experts are good, as is the case for Chamonix, in mid-December between the dashed gray lines.

During this period, the EAs should put a lot of weight on experts Q10 and Q30 which are the only good experts in this period as one can see in Figure \ref{chamonix_obs_all_experts}. However, this is not the case, since their weights remain near zero during this event as illustrated in Figure \ref{BOA_MLpol_MLprod_EWA_window1253_ech48_74056001_2020-03-30_2023-09-03_weights}. This is related to the fact that our aggregations are competing with the best expert or the best convex combination of experts and not the best sequence of experts.

This can also be explained by the fact that the weights depend on the cumulative losses. Indeed, the duration of unusual weather events - where the on average bad experts are good - are typically short compared to the number of days where these experts were bad. Thus, despite being online, the EAs need a too large number of iterations in order to compensate large cumulative losses and quickly put a lot of weight on a previously bad expert. Therefore, these EAs are not able to quickly switch a lot of weight from one expert to another.

The weights of NBM in Figure \ref{BOA_MLpol_MLprod_EWA_window1253_ech48_74056001_2020-03-30_2023-09-03_weights} are much closer to the uniform distribution than the weights of the EAs. This comes from the first step of NBM reducing the bias of the experts and making their predictions more similar and also from the constant decaying weight factor.

Hence unlike EA, which tries to ignore poor experts, NBM tries to upgrade all experts and then averages them. Therefore, in our EAs some experts are no longer used even though they are sometimes good. Whereas NBM has the opposite problem: all experts are always used but no expert is anymore adapted to predict quantiles of the temperature distribution.

\subsection{Sliding Window}
We have seen in the previous subsection that the EAs have trouble using the biased experts at the right moment because they made too often bad predictions in the past. So one might wonder whether or not it is constructive to use all the past data, especially for weather predictions, given the fact that weather shows periodic and seasonal behavior. Furthermore, focusing on recent data could help the aggregations to deal with the biased experts. Indeed, they would suffer less from their (on average) past bad performances for the weight updates.

So as in \citet{zamo_sequential_2021} and \citet{stoltz_agregation_2010} we looked at the effect of a sliding window on NBM and the EAs with the gradient trick and all experts. In order to do this, for a sliding window of $w \geqslant 1$ iterations, at each iteration $t \geqslant w+1$ , we replaced each variable of the form $X_t=\sum_{s=1}^t x_s$ by $X_t=\sum_{s=t-w}^t x_t$ and each variable of the form $Y_t= \min (y_s, 1 \leqslant s \leqslant t)$ by $Y_t= \min (y_s, t-w \leqslant s \leqslant t)$.

\begin{figure*}[h]
 \centerline{\includegraphics[width=39pc]{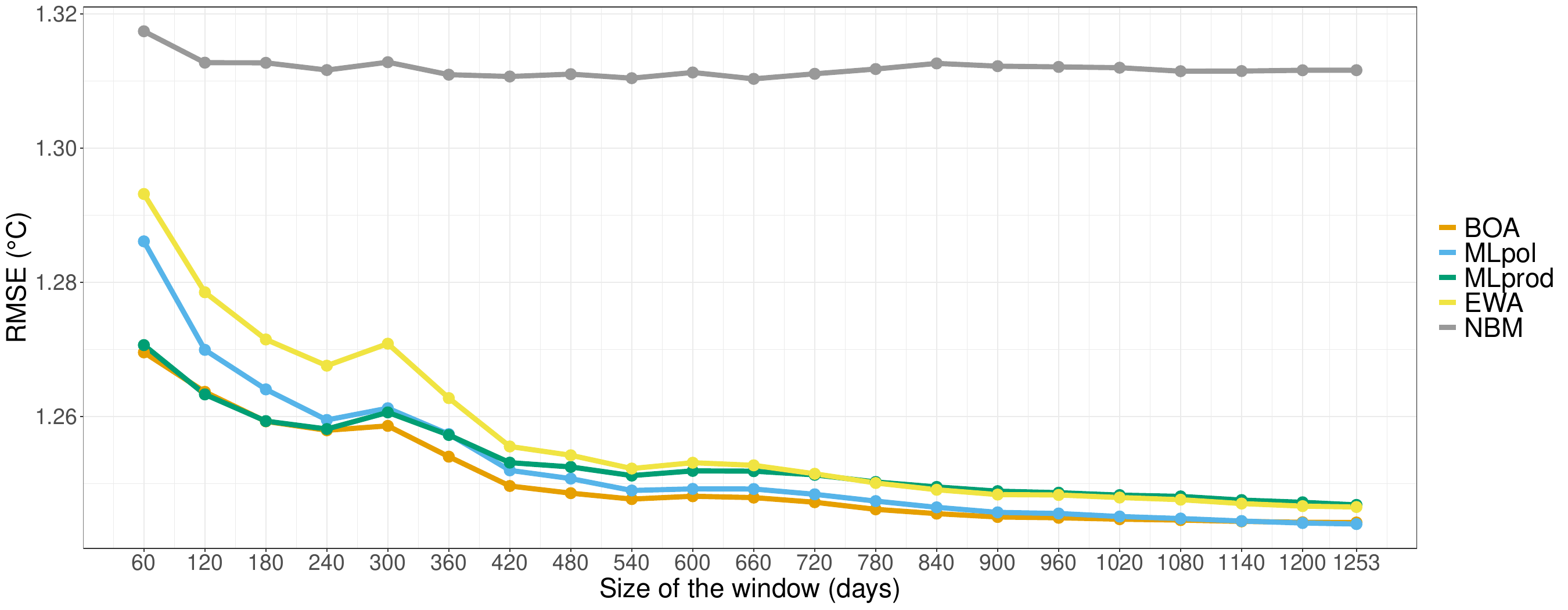}}
  \caption{RMSE over all the stations and lead times depending on the sliding windows of NBM \citep{craven_national_2020} and of the expert aggregations (with the gradient trick and all the experts): BOA \citep{wintenberger_stochastic_2024}, EWA \citep{vovk_aggregating_1990}, MLpol and MLprod \citep{gaillard_second-order_2014}.}\label{sliding_window}
\end{figure*}

In Figure \ref{sliding_window} we can observe the RMSE scores of different aggregations depending on the sliding window. It shows that in line with the regret bounds, the RMSE of the EAs improves when the sliding window increases.

The plot also shows that the EAs achieve a large reduction in RMSE with sliding windows of approximately 360 days. This is because the most relevant data for prediction – data from the same season as the day being forecasted – is included within these windows. Conversely, the RMSE increases slightly around the 300 and 600 sliding windows, where data from other seasons than the day to be predicted is incorporated.

After sliding windows of around 600 days, the decrease of the EA's RMSE is slower than previously. This suggests that using data older than one year and a half does not help the EA as much as recent data.

It is also quite impressive that even for the smallest sliding window of 60 days that we tried out, all the aggregation strategies have a better RMSE than those of the uniform aggregation (1.41°C) and the best expert in hindsight (1.37°C). This holds even for NBM, even though it performs worse than the other EAs for every sliding window tested.

The fact that the RMSE of NBM quickly stops to decrease when the sliding window increases, indicates that unlike the EAs, NBM is not able to learn from old data. This is due to the decaying weight factor which is not adaptive (unlike the learning rates of the EAs) and favors recent data \citep{craven_national_2020}.

This plot also shows that in our study, BOA and MLprod are the "fastest" learners since for the small sliding windows they have the lowest RMSE. This means that when fewer than 250 predictions are to be made, BOA or MLprod should be preferred. 

\begin{figure*} [h]
 \centerline{\includegraphics[width=39pc]{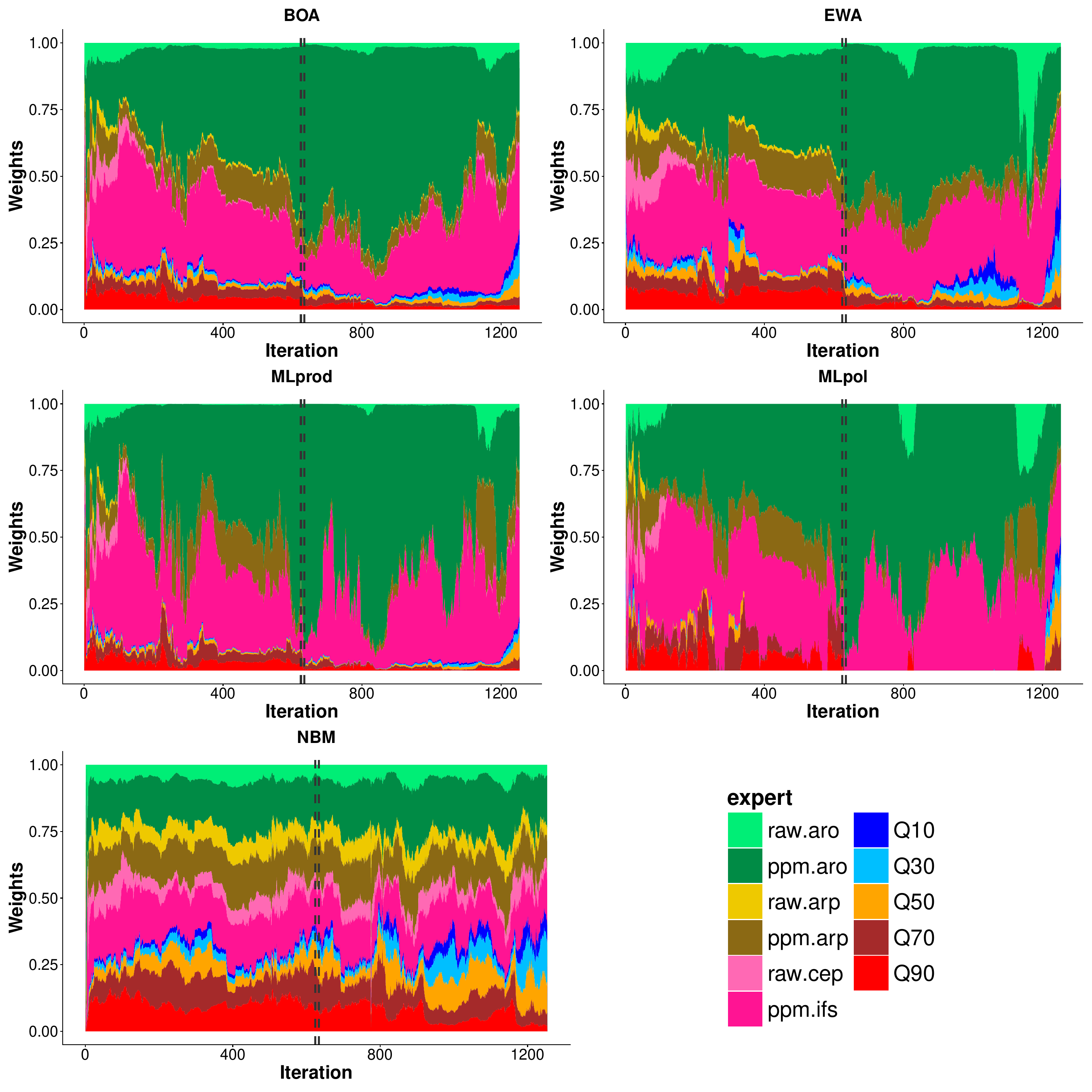}}
  \caption{Weights of the experts in Chamonix, from 2020-03-30 to 2023-09-03, for the lead time 48 hours, for different expert aggregations with the gradient trick and a sliding window of 500 days (EWA \citep{vovk_aggregating_1990}, BOA \citep{wintenberger_stochastic_2024}, MLpol and MLprod \citep{gaillard_second-order_2014}) and for NBM \citep{craven_national_2020}. The experts are the "raw" and Post Processed Models (ppm) AROME (aro), IFS (ifs), ARPEGE (arp) and the post processed quantiles of the PEARP modele (Q10, Q30, Q50, Q70, Q90). The two dashed grey - almost overlapping - lines delimit the mid-December period where all experts make too warm predictions except the Q10 and Q30 experts.}\label{BOA_MLpol_MLprod_EWA_window500_ech48_74056001_2020-03-30_2023-09-03_weights}
\end{figure*}

Finally, we can see on Figure \ref{BOA_MLpol_MLprod_EWA_window500_ech48_74056001_2020-03-30_2023-09-03_weights} that one can make the EA much more reactive by adding a sliding window. The variance of the weights is greatly increased with the sliding window compared with the weights of the EAs on Figure \ref{BOA_MLpol_MLprod_EWA_window1253_ech48_74056001_2020-03-30_2023-09-03_weights}. Even the biased experts can have large weights at any time.

But there is a trade-off between the RMSE and the reactivity through the sliding window, because reducing the sliding window certainly adds some noise.  However, the faster and larger weight variations cannot be explained solely by the added noise. 

Some of the experts have very similar predictions, hence the aggregations will quickly switch from one expert to another. This can be observed on Figure \ref{BOA_MLpol_MLprod_EWA_window500_ech48_74056001_2020-03-30_2023-09-03_weights} with the weights of the two best experts ppm.ifs and ppm.aro.

The fast variations on relatively short time scales also seem to reflect the changes in the quality of the experts over shorter time scales. Indeed, on several occasions all EAs rapidly assigned a large weight to the same expert as can be seen on Figure \ref{BOA_MLpol_MLprod_EWA_window500_ech48_74056001_2020-03-30_2023-09-03_weights}. 
For example, raw.aro gains a lot of weight just before iteration 1200 for all EAs.

However, despite greater reactivity, a sliding window of 500 days does not seem short enough to catch about 10 large consecutive errors of the on average good experts. Indeed, none of the aggregations are able to put a lot of weight on the Q10 and Q30 experts during the mid-December event in Chamonix.

By comparing Figures \ref{BOA_MLpol_MLprod_EWA_window500_ech48_74056001_2020-03-30_2023-09-03_weights} and \ref{BOA_MLpol_MLprod_EWA_window1253_ech48_74056001_2020-03-30_2023-09-03_weights}, one can see that NBM's weights with a 500 days sliding window are similar to those without a sliding window. Again, this can be explained by the decaying weight factor which does not depend on the experts nor on the iteration and the fact that the experts are similar after the bias reduction.

So, one can use a sliding window with EA to make the aggregation more reactive, whereas this seems more difficult with NBM. And overall, BOA always has the best RMSE for all sliding windows compared to the other aggregations. Hence, BOA has the fastest convergence.

\section{Conclusion}
\label{section_conclusion}
We have observed that all the EAs significantly improve the 2 meters above the ground temperature predictions when compared to all available raw and post-processed NWP models. These aggregations also outperformed the mean of all available NWP models.

It is quite remarkable that the gradient trick enables EAs with second-order regret bounds to achieve a RMSE close to that of the best fixed convex combination of experts, without knowing this oracle in advance. 

It is the power of EA to be able to improve a large set of experts - including already post-processed forecasts - by combining them. Hence, EA enables more reliable forecasting for operational purposes than the experts it relies on.

The differences between the EAs with the gradient trick are not statistically significant (except for a very small fraction of station-lead time couples in favor of BOA against EWA and MLprod). 
Nevertheless, we advise to use BOA, which has strong theoretical guarantees and the best scores in every setting that we tested.

All EAs outperformed NBM; in the case of BOA, the RMSE differences are statistically significant for about 93\% of the station-lead time pairs. And Figure \ref{boxplot_with_pearp_with_grad} showed that BOA has a better RMSE distribution for all the lead times compared to NBM.

In this study we only worked on temperature, but it would also be possible to apply EA to other variables like sea level pressure or humidity.  Moreover, while we only used raw and post-processed NWP models, one can in principle aggregate all types of models.

Therefore, it is possible to improve the predictions of several distinct models in a fully online and adaptive way with fairly simple algorithms. The computational overhead of the EA itself is negligible, approximately 1 minute on the Météo France supercomputer for roughly 2500 stations and 32 lead times. 

One can also observe that EA helps reduce retraining costs, as it is an adaptive online learning algorithm. Thereby, EA does not require to be retrained every time a model changes. If the models change and their post processed versions are not retrained immediately, the EA will adjust its weights to compensate for any degradation in the post-processed outputs. In that case, resorting to aggregation will mitigate the degradation. This is especially the case when it is used with a sliding window.

All of these EAs however struggle rapidly to really use biased experts which may be useful when the other experts are bad. To remedy this, it seems possible to find a trade-off between average scores and the reactivity of the aggregation by using a sliding window.

One may wonder however if it is possible to improve the very good scores on average of the aggregations with second order bounds whilst being more reactive. This would help to lower the $Q_{95}(|e|)$ by using the biased experts more efficiently for extreme events. 

One way could be to try to compete with the best sequence of experts such as in \citet{herbster_tracking_1998} and \citet{mourtada_efficient_2017} instead of a fixed expert as we did in this study. Another way is to use the Sleeping Expert Framework (SEF) as \citet{freund_using_1997} and \citet{devaine_forecasting_2013}. In Part II, we use the SEF to improve the EAs not only in terms of RMSE and $Q_{95}(|e|)$, but also in certain operationally critical situations.

\clearpage
\section*{Acknowledgments}
The authors would like to thank the three anonymous reviewers for their helpful comments and suggestions that greatly improved the quality of the article. The authors would like to acknowledge the support of the French Agence Nationale de la Recherche (ANR) under reference ANR20-CE40-0025-01 (T-REX project). Finally, the authors are grateful to Hermann Pfitzner for his thorough review.

\section*{Data availability}
The data used for the experiments is available at https://github.com/pfitznerl/agregation\_2025.

\section*{Appendix, Statistical significance tests}

\begin{table*}[h]
\centering
\begin{tabular}{|c|c|c|c|c|c|c|c|c|}
\hline
 & bestConvex & BOA & EWA & MLpol & MLprod & NBM & UNIFORM & best expert\\
\hline
bestConvex & & 78.1 & 79.5 & 68.0 & 74.9 & 98.7 & 98.8 & 100 \\
BOA & 0 & & 0.2 & 0 & 0.2 & 93.3 & 99.7 & 100 \\
EWA & 0 & 0 & & 0 & 0 & 87.4 & 98.7 & 100 \\
MLpol & 0 & 0 & 0 & & 0 & 88.4 & 97.6 & 100 \\
MLprod & 0 & 0 & 0 & 0 & & 81.5 & 96.1 & 100 \\
NBM & 0 & 0 & 0 & 0 & 0 & & 41.4 & 61.6 \\
UNIFORM & 0 & 0 & 0 & 0 & 0 & 0 & & 47.5 \\
bestExpert & 0 & 0 & 0 & 0 & 0 & 0.2 & 13.5 & \\
\hline
\end{tabular}
\caption{Percentage of station-lead time pairs with a statistically significative difference between the aggregation's RMSE. On row $i$ and column $j$ is the percentage of lead time/station pairs for which the one sided Diebold–Mariano test rejects the null hypothesis in favor of the aggregation of row $i$ against the aggregation of column $j$ with a Benjamini–Hochberg procedure to take account of multiple testing with a 0.05 level. The tests were run with the dm.test function of the R package forecast 8.24.0 and the p.adjust function of the R package stats 3.6.1.}
\label{test_rmse}
\end{table*}

\begin{table*}[h]
\centering
\begin{tabular}{|c|c|c|c|c|c|c|c|c|}
\hline
 & bestConvex & BOA & EWA & MLpol & MLprod & NBM & UNIFORM & best expert\\
\hline
bestConvex & & 0 & 0 & 0 & 0 & 0 & 29.5 & 58.2 \\
BOA & 0 & & 0 & 0 & 0 & 0 & 20.9 & 47.0 \\
EWA & 0 & 0 & & 0 & 0 & 0 & 19.7 & 37.9 \\
MLpol & 0 & 0 & 0 & & 0 & 0 & 19.7 & 46.4 \\
MLprod & 0 & 0 & 0 & 0 & & 0 & 18.9 & 37.4 \\
NBM & 0 & 0 & 0 & 0 & 0 & & 9.4 & 1.0 \\
UNIFORM & 0 & 0 & 0 & 0 & 0 & 0 & & 0.2 \\
bestExpert & 0 & 0 & 0 & 0 & 0 & 0 & 7.6 & \\
\hline
\end{tabular}
\caption{Percentage of station-lead time pairs with a statistically significative difference between the aggregation's $Q_{95}(|e|)$. On row $i$ and column $j$ is the percentage of station-lead time pairs for which the quantile test rejects the null hypothesis in favor of the aggregation of row $i$ against the aggregation of column $j$ with a Benjamini–Hochberg procedure to take account of multiple testing with a 0.05 level. The tests were run with the quantileTest function of the R package EnvStats 3.1.0 and the p.adjust function of the R package stats 3.6.1.}
\label{test_q95}
\end{table*}

\bibliographystyle{unsrtnat}
\bibliography{references.bib}  






\end{document}